\documentclass[12pt,a4paper,twoside]{article}
\usepackage[english]{babel}
\usepackage{amssymb}
\usepackage{inputenc}
\newcommand{\ba}{\begin{array}}
\newcommand{\ea}{\end{array}}
\usepackage{amssymb}
\begin{document}

\author{S. Albeverio $^{1},$ Sh. A. Ayupov $^{2,  *},$  K. K.
Kudaybergenov  $^3$}

\title{\bf Structure of derivations  on  various algebras of measurable operators for type I  von Neumann algebras}

\maketitle

\begin{abstract}
Given a von Neumann algebra $M$ denote by $S(M)$ and $LS(M)$ respectively
the algebras of all  measurable and  locally measurable operators  affiliated with $M.$
For a faithful normal semi-finite trace $\tau$ on $M$ let $S(M, \tau)$ (resp. $S_0(M, \tau)$) be the
algebra of all $\tau$-measurable (resp. $\tau$-compact) operators from $S(M).$ We give a complete description of all derivations on the above algebras of operators in the case of type I von Neumann algebra $M.$ In particular, we prove that if $M$ is of type I$_\infty$ then every derivation on $LS(M)$ (resp. $S(M)$ and $S(M,\tau)$)
is inner, and each derivation on $S_0(M, \tau)$ is spatial and implemented by an element from $S(M, \tau).$
\end{abstract}

\medskip
$^1$ Institut f\"{u}r Angewandte Mathematik, Universit\"{a}t Bonn,
Wegelerstr. 6, D-53115 Bonn (Germany); SFB 611, BiBoS; CERFIM
(Locarno); Acc. Arch. (USI), \emph{albeverio@uni-bonn.de}

$^2$ Institute of Mathematics and information  technologies,
Uzbekistan Academy of Sciences, F. Hodjaev str. 29, 100125,
Tashkent (Uzbekistan), e-mail: \emph{sh\_ayupov@mail.ru}

 $^{3}$ Karakalpak state university, Ch. Abdirov str. 1, 742012, Nukus (Uzbekistan),
e-mail: \emph{karim2006@mail.ru}

 \medskip \textbf{AMS Subject Classifications (2000):} 46L57, 46L50, 46L55,
46L60

\textbf{Key words:}  von Neumann algebras,  non commutative
integration,  measurable operator, locally measurable operator, $\tau$-measurable operator, $\tau$-compact operator, type I
von Neumann algebra, derivation, spatial derivation.

* Corresponding author
\newpage

\begin{center}
{\bf Introduction}
\end{center}

Derivations on unbounded operator algebras, in particular on
various algebras of measurable operators affiliated with von
Neumann algebras, appear to be a very attractive special case of
the general theory of unbounded derivations on operator algebras.
The present paper continues the series of papers of the authors
\cite{Alb1}-\cite{Alb3} devoted to the study and a description of derivations on the
algebra $LS(M)$ of locally measurable operators with respect to a
von Neumann algebra $M$ and on various subalgebras of $LS(M).$

Let  $A$ be an algebra over the complex number. A linear operator
$D:A\rightarrow A$ is called a \emph{derivation} if it satisfies the
identity  $D(xy)=D(x)y+xD(y)$ for all  $x, y\in A$ (Leibniz rule).
Each element  $a\in A$ defines a derivation  $D_a$ on $A$ given as
$D_a(x)=ax-xa,\,x\in A.$ Such derivations $D_a$ are said to be
\emph{inner derivations}. If the element  $a$ implementing the
derivation   $D_a$ on $A,$ belongs to a larger algebra  $B,$
containing  $A$ (as a proper ideal as usual) then $D_a$ is called
a  \emph{spatial derivation}.

In the particular case where $A$ is commutative, inner derivations
are identically zero, i.e. trivial. One of the main problems in
the theory of derivations is automatic innerness or spatialness of
derivations and the existence of non inner derivations (in
particular, non trivial derivations on commutative algebras).

On this way  A.~F.~Ber, F.~A.~Sukochev, V.~I.~Chilin~\cite{Ber}
obtained necessary and sufficient conditions for the existence of
non trivial derivations on commutative regular algebras. In
particular they have proved that the algebra  $L^{0}(0, 1)$ of all
(classes of equivalence of) complex measurable functions on  the  interval $(0,
1)$
admits non trivial derivations. Independently
A.~G.~Kusraev~\cite{Kus1} by means of Boolean-valued analysis has
established necessary and sufficient conditions for the existence
of non trivial derivations and automorphisms on universally complete complex
$f$-algebras. In particular he has also proved the existence of
non trivial derivations and automorphisms on  $L^{0}(0, 1).$ It is clear that these derivations
are discontinuous in the measure topology, and therefore they are neither
inner nor spatial. It seems that the existence of such pathological example of derivations deeply
depends on the commutativity of the underlying von Neumann algebra $M.$
In this connection the present authors have initiated the study of the above problems
in the non commutative case \cite{Alb1}-\cite{Alb4}, by considering derivations
on the algebra $LS(M)$ of all locally measurable operators with
respect to a semi-finite von Neumann algebra $M$ and on various
subalgebras of  $LS(M).$ Recently another  approach to similar problems in the framework of type I $AW^{*}$-algebras has been outlined in \cite{Gut}.

The main result of the paper \cite{Alb1} states that if
$M$ is a type I von Neumann algebra, then every derivation $D$ on $LS(M)$
which is identically zero on the center $Z$ of the von Neumann algebra
$M$ (i.e. which is $Z$-linear) is automatically inner, i.e. $D(x)=ax-xa$ for an appropriate
$a\in LS(M).$ In  \cite[Example 3.8]{Alb1} we also gave a construction
of non inner derivations $D_\delta$ on the algebra $LS(M)$ for type I$_{fin}$ von Neumann
algebra $M$ with non atomic center $Z$, where $\delta$ is a non
trivial derivation on the algebra $LS(Z)$ (i.e. on the center of $LS(M)$) which is isomorphic with the algebra $L^{0}(\Omega, \Sigma, \mu)$
of all measurable functions on a non atomic measure space $(\Omega, \Sigma, \mu).$

The main idea of the mentioned construction is the following.

Let $A$ be a commutative algebra and let $M_n(A)$ be the algebra of $n\times n$ matrices over $A.$
If  $e_{i,j},\,i,j=\overline{1, n},$
are the matrix units in  $M_n(A),$ then each element
$x\in M_n(A)$ has the form
 $$x=\sum\limits_{i,j=1}^{n}\lambda_{i,j}e_{i,j},\,\lambda_{i,j}\in A,\,i,j=\overline{1, n}.$$
Let  $\delta:A\rightarrow A$ be a
derivation. Setting

 \begin{equation}
  D_{\delta}(\sum\limits_{i,j=1}^{n}\lambda_{i,j}e_{i,j})=
 \sum\limits_{i,j=1}^{n}\delta(\lambda_{i,j})e_{i,j}
\end{equation}
 we obtain a well-defined linear operator
 $D_\delta$ on the algebra $M_n(A).$ Moreover
 $D_\delta$ is a derivation on the algebra  $M_n(A)$
 and its restriction onto the center of the algebra  $M_n(A)$ coincides with the given $\delta.$

 In papers \cite{Alb2},  \cite{Alb4} we considered similar problems for derivations on the algebra $S_0(M, \tau)$
of $\tau$-compact operators with respect to a type I von Neumann algebra $M$ with a faithful normal semi-finite trace $\tau,$ and obtained necessary and sufficient conditions for derivations to be spatial.
 In  \cite{Alb3} we have proved spatialness
of all derivations on the non commutative Arens algebra
$L^{\omega}(M, \tau)$ associated with an arbitrary von Neumann
algebra  $M$ and a faithful normal semi-finite trace  $\tau.$ Moreover if the trace $\tau$
is finite then every derivation on $L^{\omega}(M, \tau)$ is inner.

In the present paper we give a complete description of all derivations on the algebra $LS(M)$
of all locally measurable operators affiliated with a type I von Neumann algebra $M$, and also
on its subalgebras $S(M)$ -- of measurable operators, $S(M, \tau)$ of $\tau$-measurable operators
and on $S_0(M, \tau)$ of all $\tau$-compact operators with respect to $M$,
where $\tau$ is a faithful normal semi-finite trace on $M.$
We prove that the above mentioned construction of derivations $D_\delta$ from \cite{Alb1} gives the general form of pathological derivations on these algebras and these exist only in the type
I$_{fin}$ case, while for type I$_\infty$ von Neumann algebras $M$ all derivations on $LS(M),$ $S(M)$
and $S(M, \tau)$ are inner and for $S_0(M, \tau)$ they  are spatial. Moreover we prove that an
arbitrary derivation $D$ on each of these algebras can be uniquely decomposed into the sum $D=D_a+D_\delta$
where the derivation $D_a$ is inner (for $LS(M),$ $S(M)$
and $S(M, \tau)$) or spatial (for $S_0(M, \tau)$) while the derivation $D_\delta$ is constructed in the above mentioned manner from a non trivial  derivation $\delta$ on the center of the corresponding algebra.

In section 1 we give necessary definition and preliminaries from the theory of measurable operators and Hilbert -- Kaplansky modules.

In section 2 we describe derivations on the algebra $LS(M)$ of all locally measurable
operators for a type I von Neumann algebra $M.$

Sections 3 and 4 are devoted to derivation respectively on the algebra $S(M)$ of all measurable operators and on the algebra $S(M, \tau)$
of all $\tau$-measurable operators with respect to $M,$ where $M$ is a type I von Neumann algebra and $\tau$
is a faithful normal semi-finite trace on $M.$

 In Section 5 we give the solution of the problem for derivations on the algebra $S_0(M, \tau)$
of all $\tau$-compact operators affiliated with a type I von Neumann algebra $M$ and a faithful normal semi-finite trace $\tau.$

Finally, section 6 contains an application of the above results to the description of the first cohomology group for the considered algebras.
\begin{center}
{\bf 1. Preliminaries}
\end{center}

Let  $H$ be a complex Hilbert space and let  $B(H)$ be the algebra
of all bounded linear operators on   $H.$ Consider a von Neumann
algebra $M$  in $B(H)$ with the operator norm $\|\cdot\|_M.$ Denote by  $P(M)$ the lattice of projections in $M.$

A linear subspace  $\mathcal{D}$ in  $H$ is said to be
\emph{affiliated} with  $M$ (denoted as  $\mathcal{D}\eta M$), if
$u(\mathcal{D})\subset \mathcal{D}$ for every unitary  $u$ from
the commutant
$$M'=\{y\in B(H):xy=yx, \,\forall x\in M\}$$ of the von Neumann algebra $M.$

A linear operator  $x$ on  $H$ with the domain  $\mathcal{D}(x)$
is said to be \emph{affiliated} with  $M$ (denoted as  $x\eta M$) if
$\mathcal{D}(x)\eta M$ and $u(x(\xi))=x(u(\xi))$
 for all  $\xi\in
\mathcal{D}(x).$

A linear subspace $\mathcal{D}$ in $H$ is said to be \emph{strongly
dense} in  $H$ with respect to the von Neumann algebra  $M,$ if

1) $\mathcal{D}\eta M;$

2) there exists a sequence of projections
$\{p_n\}_{n=1}^{\infty}$ in $P(M)$  such that
$p_n\uparrow\textbf{1},$ $p_n(H)\subset \mathcal{D}$ and
$p^{\perp}_n=\textbf{1}-p_n$ is finite in  $M$ for all
$n\in\mathbb{N},$ where $\textbf{1}$ is the identity in $M.$

A closed linear operator  $x$ acting in the Hilbert space $H$ is said to be
\emph{measurable} with respect to the von Neumann algebra  $M,$ if
 $x\eta M$ and $\mathcal{D}(x)$ is strongly dense in  $H.$ Denote by
 $S(M)$ the set of all measurable operators with respect to
 $M.$

A closed linear operator $x$ in  $H$  is said to
be \emph{locally measurable} with respect to the von Neumann
algebra $M,$ if $x\eta M$ and there exists a sequence
$\{z_n\}_{n=1}^{\infty}$ of central projections in $M$ such that
$z_n\uparrow\textbf{1}$ and $z_nx \in S(M)$ for all
$n\in\mathbb{N}.$

It is well-known \cite{Mur1} that the set $LS(M)$ of all locally
measurable operators with respect to $M$ is a unital *-algebra when equipped with
the algebraic operations of strong addition and multiplication and taking the adjoint of an operator.

 Let   $\tau$ be a faithful normal semi-finite trace on $M.$ We recall that a closed linear operator
  $x$ is said to be  $\tau$\emph{-measurable} with respect to the von Neumann algebra
   $M,$ if  $x\eta M$ and   $\mathcal{D}(x)$ is
  $\tau$-dense in  $H,$ i.e. $\mathcal{D}(x)\eta M$ and given   $\varepsilon>0$
  there exists a projection   $p\in M$ such that   $p(H)\subset\mathcal{D}(x)$ and $\tau(p^{\perp})<\varepsilon.$
   The set $S(M,\tau)$ of all   $\tau$-measurable operators with respect to  $M$
    is a solid  *-subalgebra in $S(M)$  (see \cite{Nel}).

    Consider the topology  $t_{\tau}$ of convergence in measure or \emph{measure topology}
    on $S(M, \tau),$ which is defined by
 the following neighborhoods of zero:
$$V(\varepsilon, \delta)=\{x\in S(M, \tau): \exists e\in P(M), \tau(e^{\perp})\leq\delta, xe\in
M,  \|xe\|_{M}\leq\varepsilon\},$$  where $\varepsilon, \delta$ are positive numbers.

 It is well-known
\cite{Nel} that $S(M, \tau)$ equipped with the measure topology is a
complete metrizable topological *-algebra.

An element  $x$ of the algebra  $S(M, \tau)$ is said to be  $\tau$-\emph{compact}, if given any
 $\varepsilon >0$ there exists a projection $p\in P(M)$
such that  $\tau(p^{\perp})<\infty,\,xp\in M$ and
$\|xp\|_M<\varepsilon.$ The set $S_0(M, \tau)$ of all $\tau$-compact operators
is an  *-ideal in the algebra
$S(M, \tau)$
(see \cite{Mur1}).

It should be noted that the algebra  of $\tau$-compact operators were considered by Yeadon \cite{Yea} and  Fack and Kosaki \cite{Fac} and one of
the original definitions was the following: an operator $x\in S(M, \tau)$ is said to be $\tau$-compact if
$$\lim\limits_{t\rightarrow\infty}\mu_{t}(x)=0,$$
where $\mu_{t}(x)=\inf\{\lambda>0: \tau(e_{\lambda}^{\perp})\leq t\}$ and
$\{e_{\lambda}\}_{\lambda>0}$ is the spectral resolution of $|x|.$  The equivalence of this definition to the one given above was proved in
 \cite{Str}.

Note that if the trace $\tau$ is a finite then
$$S_0(M, \tau)=S(M, \tau)=S(M)=LS(M).$$

The following result describes one of the most important
properties of the algebra  $LS(M)$ (see \cite{Mur1},
\cite{Sai}).

\textbf{Proposition 1.1.} \emph{Suppose that the von Neumann
algebra  $M$ is the  $C^{\ast}$-product of the von Neumann
algebras  $M_i,$ $i\in I,$ where $I$ is an arbitrary set of
indices, i.e.}
$$M=\bigoplus\limits_{i\in I}M_i=
\{\{x_i\}_{i\in I}:x_i\in M_i, i\in I, \sup\limits_{i\in
I}\|x_i\|_{M_i}<\infty\}$$ \emph{with coordinate-wise algebraic
operations and involution and with the $C^{\ast}$-norm
$\|\{x_i\}_{i\in I}\|_{ M}=\sup\limits_{i\in I}\|x_i\|_{M_i}.$
Then the algebra  $LS(M)$ is *-isomorphic to the algebra
$\prod\limits_{i\in I}LS(M_i)$ (with the coordinate-wise
operations and involution), i.e.}
$$LS(M)\cong\prod\limits_{i\in I}LS(M_i)$$
($\cong$ denoting *-isomorphism of algebras).

It should be noted that similar isomorphisms are not valid in general for the algebras
 $S(M),$ $S(M, \tau)$ and  $S_0(M, \tau)$  (see \cite{Mur1}).

Proposition 1.1 implies that given any family  $\{z_i\}_{i\in I}$
of mutually orthogonal central projections in $M$ with
$\bigvee\limits_{i\in I}z_i=\textbf{1}$ and a  family
of elements $\{x_i\}_{i\in I}$ in $LS(M)$ there exists a unique element $x\in LS(M)$
such that $z_i x=z_i x_i$ for all $i\in I.$ This element is
denoted by  $x=\sum\limits_{i\in I}z_i x_i.$

It is well-known \cite{Seg} that every commutative von Neumann algebra
 $M$
is *-isomorphic to the algebra   $L^{\infty}(\Omega)=L^{\infty}(\Omega, \Sigma,
\mu)$ of all (classes of equivalence of) complex essentially bounded measurable functions on a measure space
$(\Omega, \Sigma, \mu)$ and in this case  $LS(M)=S(M)\cong
L^{0}(\Omega),$   where $L^{0}(\Omega)=L^{0}(\Omega, \Sigma,
\mu)$ the algebra of all (classes of equivalence of) complex measurable functions on
$(\Omega, \Sigma, \mu).$

Further we shall need the following remarkable description of centers of the algebras $S(M),$ $S(M, \tau)$ and $S_0(M, \tau)$ for type I$_\infty$ von Neumann algebras.

\textbf{Proposition 1.2.} \emph{Let  $M$ be a type  $I_{\infty}$ von
Neumann algebra with the center $Z.$ Then}

    \emph{a) the centers of the algebras $S(M)$ and  $S(M,\tau)$ coincide with $Z;$}

    \emph{b) the center of the algebra  $S_{0}(M,\tau)$ is trivial, i.e.  $Z(S_{0}(M,\tau))=\{0\}.$}

   Proof. a) Suppose that  $z\in S(M), \ z\geq 0,$  is a central element and let
     $z=\int\limits_{0}^{\infty}\lambda\,de_{\lambda}$ be its spectral resolution. Then  $e_{\lambda}\in Z$ for all
     $\lambda>0.$ Assume that  $e_{n}^{\bot}\neq 0$ for all  $n\in\mathbb{N}.$
     Since $M$ is of type  I$_{\infty},$ $Z$ does not contain non-zero finite projections. Thus
   $e_{n}^{\bot}$ is infinite for all  $n\in\mathbb{N},$ which
   contradicts  the condition  $z\in S(M).$ Therefore there
   exists  $n_{0}\in\mathbb{N}$ such that  $e_{n}^{\bot}=0$ for all  $n\geq n_0,$
  i.e. $z\leq n_0\textbf{1}.$ This means that   $z\in Z,$ i.e. $Z(S(M))=Z.$ Similarly $Z(S(M, \tau))=Z.$

  b) Let $z\in Z(S_0(M, \tau)), \ z\geq 0.$ Take a projection $p\in
  M$ with $\tau(p)<\infty.$ Then $p\in S_0(M, \tau)$ and therefore  $zp=pz.$
Since  $M$ is semi-finite this implies that $zp=pz$ for all $p\in
P(M).$ Since the linear span of  $P(M)$ is dense in $S(M, \tau)$
in the measure topology, we have that  $zx=xz$ for all $x\in S(M,
\tau),$ i.e. $z\in Z(S(M, \tau))=Z.$

    Suppose that
    $z=\int\limits_{0}^{\infty}\lambda\,de_{\lambda}$ is the
    spectral resolution of $z.$ Then
     $e_{\lambda}\in Z$ for all  $\lambda>0.$ Since
$z\in S_0(M, \tau)$ we have that $e_{\lambda}^{\perp}$ is a
finite projection for all $\lambda>0.$ But $M$ does not contain
any non zero central finite projection, because it is of type
 I$_{\infty}.$ Therefore  $e_{\lambda}^{\perp}=0$ for all
 $\lambda>0,$ i.e.  $z=0.$ Thus  $Z(S_0(M,
    \tau))=\{0\}.$ The proof is complete. $\blacksquare$

Now let us recall some notions and results from the theory of
Hilbert -- Kaplansky modules (for details we refer to \cite{Kap}, \cite{Kus2}).

Let $(\Omega, \Sigma, \mu)$ be a measure space and let  $H$ be a Hilbert space.  A map  $s:\Omega\rightarrow H$ is said
to be simple, if
$s(\omega)=\sum\limits_{k=1}^{n}\chi_{A_{k}}(\omega)c_k,$ where
$A_k\in\Sigma, A_i\cap A_j=\emptyset, \,i\neq j,\, \,c_k\in
H,\,k=\overline{1, n},\, n\in\mathbb{N}.$ A map $u:\Omega\rightarrow
H$ is said to be measurable, if there is a sequence  $(s_n)$ of
simple maps such that $\|s_n(\omega)-u(\omega)\|\rightarrow0$
almost everywhere on any $A\in\sum$ with $\mu(A)<\infty.$

Let $\mathcal{L}(\Omega, H)$ be the set of all measurable maps from
$\Omega$ into $H,$ and let $L^{0}(\Omega, H)$ denote the space of
all equivalence classes  with respect to the equality almost
everywhere. Denote by $\hat{u}$ the equivalence class from
$L^{0}(\Omega, H)$ which contains the measurable map $u\in
\mathcal{L}(\Omega, H).$ Further we shall identify the element $u\in
\mathcal{L}(\Omega, H)$ and the class $\hat{u}.$ Note that the
function  $\omega \rightarrow \|u(\omega)\|$
       is measurable for any $u\in \mathcal{L}(\Omega, H).$ The equivalence class containing the function
              $\|u(\omega)\|$ is denoted by
       $\|\hat{u}\|$. For  $\hat{u}, \hat{v}\in L^{0}(\Omega, H), \lambda\in L^{0}(\Omega)$ put
$\hat{u}+\hat{v}=\widehat{u(\omega)+v(\omega)},
\lambda\hat{u}=\widehat{\lambda(\omega) u(\omega)}.$
Equipped with the  $L^{0}(\Omega)$-valued inner product
   $$\langle x, y \rangle=\langle x(\omega), y(\omega) \rangle_{H},$$ where  $\langle \cdot, \cdot\rangle_{H}$
   in the inner product in  $H,$
      $L^{0}(\Omega, H)$ becomes a Hilbert~-- Kaplansky module over $L^{0}(\Omega).$
       The space   $$L^{\infty}(\Omega, H)=\{x\in L^{0}(\Omega, H): \langle x,x\rangle\in L^{\infty}(\Omega)\}$$
      is a Hilbert~-- Kaplansky module over  $L^{\infty}(\Omega).$ Denote by   $B(L^{0}(\Omega, H))$ the
      algebra of all
              $L^{0}(\Omega)$-bounded
       $L^{0}(\Omega)$-linear operators on  $L^{0}(\Omega, H)$ and denote by
        $B(L^{\infty}(\Omega, H))$ the algebra of all
        $L^{\infty}(\Omega)$-bounded $L^{\infty}(\Omega)$-linear operators on
        $L^{\infty}(\Omega, H).$

Now consider a von Neumann algebra  $M$ which is homogeneous of type
I$_{\alpha}$ with the center  $L^{\infty}(\Omega),$ where $\alpha$ is a cardinal number.
Then    $M$ is *-isomorphic to the algebra
  $B(L^{\infty}(\Omega, H)),$ where $\dim H=\alpha,$ while the algebra  $LS(M)$
  is *-isomorphic to $B(L^{0}(\Omega, H))$ (see for details \cite{Alb1}).

It is known  \cite{Tak} that given a type I von Neumann algebra
$M$ there exists a unique (cardinal-indexed) family  of central
orthogonal projections  $(q_{\alpha})_{\alpha\in J}$ in
$P(M)$ with $\sum\limits_{\alpha\in
J}q_{\alpha}=\textbf{1}$ such that $q_{\alpha}M$ is a homogeneous
type I$_\alpha$ von Neumann algebra, i.e.
 $q_\alpha
M\cong B(L^{\infty}(\Omega_{\alpha}, H_{\alpha}))$ with $\dim
H_\alpha=\alpha$  and
$$M\cong\bigoplus\limits_{\alpha\in J}B(L^{\infty}(\Omega_{\alpha}, H_{\alpha})).$$
The direct product
$$\prod\limits_{\alpha\in J}L^{0}(\Omega_{\alpha}, H_{\alpha})$$
equipped with the coordinate-wise algebraic operations and inner product forms
a Hilbert~--- Kaplansky module over $L^{0}(\Omega)
\cong\prod\limits_{\alpha\in J}L^{0}(\Omega_{\alpha}).$

In  \cite{Alb1} we have proved that if the von Neumann algebra
$M$ is *-isomorphic with $\bigoplus\limits_{\alpha\in J}B(L^{\infty}(\Omega_{\alpha},
H_{\alpha}))$ then the algebra   $LS(M)$ is  *-isomorphic with   $B(\prod\limits_{\alpha\in J}L^{0}(\Omega_{\alpha},
H_{\alpha})).$ Therefore there exists a map $||\cdot||:LS(M)\rightarrow L^{0}(\Omega)$ such that for all  $x,y\in
LS(M), \lambda\in L^{0}(\Omega)$ one has
$$||x||\geq0, ||x||=0 \Leftrightarrow x=0;$$
$$||\lambda x||=|\lambda|||x||;$$
$$||x+y||\leq||x||+||y||;$$
$$||xy||\leq||x||||y||;$$
$$||xx^{*}||=||x||^{2}.$$
This map $||\cdot||:LS(M)\rightarrow L^{0}(\Omega)$ is called the \emph{center-valued} norm on $LS(M).$

\begin{center} \textbf{2. Derivations on the algebra $LS(M)$ }\end{center}

In this section we shall give a complete description of derivations on the algebra $LS(M)$ of all locally measurable operators affiliated with a type I
von Neumann algebra $M.$ It is clear that if a derivation $D$ on $LS(M)$ is inner then it is $Z$-linear, i.e. $D(\lambda x)=\lambda D(x)$
for all $\lambda\in Z,$ $x\in LS(M),$ where $Z$ is the center of the von Neumann algebra $M.$ The following main result of \cite{Alb1} asserts that the
converse is also true.

\textbf{Theorem 2.1.} \emph{Let  $M$ be a type $I$ von Neumann algebra with the center $Z.$
 Then every $Z$-linear derivation $D$ on the algebra $LS(M)$ is inner.}

Proof. (see \cite[Theorem 3.2]{Alb1}). $\blacksquare$

We are now in position to consider arbitrary (non $Z$-linear, in general)
derivations on $LS(M)$. The following simple but important remark is crucial in our further considerations.

\textbf{Remark 1.}
Let  $A$ be an algebra with the center  $Z$ and let  $D:A\rightarrow
A$ be a derivation. Given any  $x\in A$ and a central element $\lambda \in
Z$  we have
$$D(\lambda x)=D(\lambda)x+\lambda D(x)$$
and
$$D(x\lambda)=D(x)\lambda +xD(\lambda).$$
Since  $\lambda x=x\lambda $ and $\lambda D(x)=D(x)\lambda,$ it follows that  $D(\lambda)x=xD(\lambda)$ for
any $\lambda\in A.$ This means that  $D(\lambda)\in Z,$ i.e. $D(Z)\subseteq Z.$
Therefore  given any derivation $D$ on the algebra $A$ we can
consider its restriction $\delta:Z\rightarrow Z.$

Now let  $M$ be a homogeneous von Neumann algebra of type $I_{n},
n \in \mathbb{N}$,  with the center $Z.$ Then  the algebra $M$ is *-isomorphic
with the algebra $M_n(Z)$ of all  $n\times n$- matrices over $Z,$
and the algebra  $LS(M)=S(M)$ is *-isomorphic with the algebra
 $M_n(S(Z))$ of all  $n\times n$ matrices over
$S(Z),$ where $S(Z)$ is the algebra of measurable operators for the commutative von Neumann algebra $Z$.

The algebra $LS(Z)=S(Z)$ is isomorphic to the algebra $L^{0}(\Omega)=L^{}(\Omega, \Sigma, \mu)$
 of all measurable functions on a measure space
 (see section 2) and therefore it admits (in non atomic cases) non zero derivations
  (see \cite{Ber},  \cite{Kus1}).

Let  $\delta:S(Z)\rightarrow S(Z)$ be a
derivation and  $D_\delta$ be a derivation on the algebra $M_n(S(Z))$ defined by (1) in Introduction.

The following lemma describes the structure of an arbitrary derivation on the algebra of locally measurable
operators for homogeneous type  I$_n,$ $n\in\mathbb{N},$ von
Neumann algebras.

 \textbf{Lemma 2.2.} \emph{Let  $M$ be a homogenous von Neumann algebra of type
  $I_{n}, n \in \mathbb{N}.$
Every derivation  $D$ on the algebra $LS(M)$ can be uniquely represented as a sum
  $$D=D_{a}+D_{\delta ,}$$ where  $D_{a}$ is an inner derivation implemented by an element  $a\in LS(M)$
while $D_{\delta} $ is the derivation of the form (1) generated by
a derivation $\delta$ on the center of $LS(M)$ identified
with $S(Z)$.}

Proof. Let  $D$ be an arbitrary derivation on the algebra $LS(M)\cong M_n(S(Z)).$ Consider its restriction $\delta$
onto the center  $S(Z)$ of this algebra, and let
$D_\delta$ be the derivation on the algebra $M_n(S(Z))$
constructed as in (1). Put $D_1=D-D_\delta.$ Given any  $\lambda\in
S(Z)$ we have
$$D_1(\lambda)=D(\lambda)-D_\delta(\lambda)=D(\lambda)-D(\lambda)=0,$$
i.e. $D_1$ is identically zero on  $S(Z).$ Therefore $D_1$ is $Z$-linear and by Theorem 2.1 we obtain that  $D_1$
is  inner derivation and thus $D_1=D_a$ for an appropriate $a\in
M_n(S(Z)).$ Therefore $D=D_a+D_\delta.$

Suppose that
$$D=D_{a_{1}}+D_{\delta_{1}}=D_{a_{2}}+D_{\delta_{2}}.$$ Then
$D_{a_{1}}-D_{a_{2}}=D_{\delta_{2}}-D_{\delta_{1}}.$ Since
$D_{a_{1}}-D_{a_{2}}$ is identically zero on the center of the
algebra $M_n(S(Z))$ this implies that
$D_{\delta_{2}}-D_{\delta_{1}}$ is also identically zero on the
center of  $M_n(S(Z)).$ This means that
$\delta_{1}=\delta_{2},$ and therefore  $D_{a_{1}}=D_{a_{2}},$
i.e. the decomposition of $D$ is unique. The proof is complete. $\blacksquare$

Now let  $M$ be an arbitrary finite von Neumann algebra of type I
with the center $Z.$ There exists a family  $\{z_n\}_{n\in F},$
$F\subseteq\mathbb{N},$ of central projections from $M$ with
$\sup\limits_{n\in F}z_n=\textbf{1}$ such that the algebra  $M$ is
*-isomorphic with the  $C^{*}$-product of von Neumann algebras
$z_n M$ of type  I$_{n}$ respectively, $n\in F,$ i.e.
$$M\cong\bigoplus\limits_{n\in F}z_n M.$$
By Proposition 1.1 we have that
$$LS(M)\cong\prod\limits_{n\in F}LS(z_n M).$$

Suppose that    $D$ is a derivation on  $LS(M),$ and
$\delta$ is its restriction onto its center  $S(Z).$
Since  $\delta$ maps each  $z_nS(Z)\cong Z(LS(z_n M))$ into itself,  $\delta$ generates a derivation $\delta_n$
on $z_nS(Z)$ for each $n\in F.$

Let     $D_{\delta_n}$ be the derivation on the matrix algebra
$M_n(z_nZ(LS(M)))\cong LS(z_nM)$ defined as in
(1). Put
\begin{equation}
D_\delta(\{x_n\}_{n\in F})=\{D_{\delta_n}(x_n)\},\,\{x_n\}_{n\in F}\in LS(M).
\end{equation}
Then the map  $D$  is a derivation on $LS(M).$

Now Lemma 2.2 implies the following result:

\textbf{Lemma 2.3.} \emph{Let  $M$ be a finite von Neumann algebra
of type I. Each
derivation  $D$ on the algebra  $LS(M)$ can be uniquely
represented in the form
$$D=D_{a}+D_{\delta ,}$$
where  $D_{a}$ is an inner derivation implemented by an element
$a\in LS(M),$ and $D_{\delta} $ is a derivation given as
(2).}

In order to consider the case of type I$_{\infty}$ von Neumann
algebra we need
some auxiliary results concerning derivations on the algebra $L^{0}(\Omega)=L^{}(\Omega, \Sigma, \mu).$

Recall that a net $\{\lambda_\alpha\}$ in $L^{0}(\Omega)$ $(o)$-converges to
$\lambda\in L^{0}(\Omega)$ if there exists a net $\{\xi_\alpha\}$ monotone
decreasing to zero such that
$|\lambda_\alpha-\lambda|\leq\xi_\alpha$ for all $\alpha.$

Denote by $\nabla$ the complete Boolean algebra of all idempotents
from $L^{0}(\Omega),$ i. e. $\nabla=\{\tilde{\chi}_{A}: A\in\Sigma\},$ where
$\tilde{\chi}_{A}$ is  the element from $L^{0}(\Omega)$ which contains the
characteristic function of the set $A.$ A \emph{partition of the
unit} in $\nabla$ is a family $(\pi_\alpha)$ of orthogonal
idempotents from $\nabla$ such that
$\bigvee\limits_{\alpha}\pi_\alpha=\textbf{1}.$

 \textbf{Lemma 2.4.} \emph{Any derivation  $\delta$ on the algebra
 $L^{0}(\Omega)$ commutes with the mixing operation on $L^{0}(\Omega),$ i.e.
$$\delta(\sum\limits_{\alpha}\pi_{\alpha}\lambda_{\alpha})=
\pi_{\alpha}\delta(\sum\limits_{\alpha}\lambda_{\alpha})$$ for an arbitrary family
$(\lambda_\alpha)\subset L^{0}(\Omega)$ and   any partition $\{\pi_{\alpha}\}$ of
the unit in $\nabla.$}

Proof. Consider a family  $\{\lambda_\alpha\}$ in $L^{0}(\Omega)$ and a
partition of the unit  $\{\pi_{\alpha}\}$ in $\nabla\subset L^{0}(\Omega).$
Since  $\delta(\pi)=0$ for any idempotent  $\pi\in\nabla,$ we have
$\delta(\pi_\alpha)=0$ for all $\alpha$ and thus
$\delta(\pi_{\alpha}\lambda)=\pi_{\alpha}\delta(\lambda)$ for any
$\lambda\in L^{0}(\Omega).$ Therefore for each $\pi_{\alpha_{0}}$ from the
given partition of the unit we have
$$\pi_{\alpha_{0}}\delta(\sum\limits_{\alpha}\pi_{\alpha}\lambda_{\alpha})=
\delta(\pi_{\alpha_{0}}\sum\limits_{\alpha}\pi_{\alpha}\lambda_{\alpha})=
\delta(\pi_{\alpha_{0}}\lambda_{\alpha_{0}})=
\pi_{\alpha_{0}}\delta(\lambda_{\alpha_0}).$$ By taking the sum over
all $\alpha_0$
 we obtain $$\delta(\sum\limits_{\alpha}\pi_{\alpha}\lambda_{\alpha})=
\sum\limits_{\alpha}\pi_{\alpha}\delta(\lambda_{\alpha}).$$ The
proof is complete. $\blacksquare$

 Recall \cite{Kus2} that a subset $K\subset L^{0}(\Omega)$ is called \emph{cyclic}, if
       $\sum\limits_{\alpha\in J}\pi_{\alpha}u_{\alpha}\in K$ for each family $(u_{\alpha})_{\alpha\in J}\subset K$ and
       for any partition  of the unit $(\pi_{\alpha})_{\alpha\in J}$ in $\nabla.$

\textbf{Lemma 2.5.} \emph{Given any non trivial derivation
$\delta:L^{0}(\Omega)\rightarrow L^{0}(\Omega)$ there exist a sequence
$\{\lambda_n\}_{n=1}^{\infty}$ in $L^{\infty}(\Omega)$ with} $|\lambda_n|\leq
\textbf{1},\,n\in \mathbb{N},$ \emph{and an idempotent $\pi\in
\nabla,\,\pi\neq 0$ such that
$$|\delta(\lambda_n)|\geq n\pi$$
for all $n\in \mathbb{N}.$}

Proof. Suppose that the set  $\{\delta(\lambda): \lambda\in L^{0}(\Omega),
|\lambda|\leq\textbf{1}\}$ is order bounded in  $L^{0}(\Omega).$  Then
$\delta$ maps any uniformly convergent sequence in
$L^{\infty}(\Omega)$ to an $(o)$-convergent sequence in $L^{0}(\Omega).$
The algebra $L^{\infty}(\Omega)$ coincides with the uniform
closure of the linear span of idempotents from $\nabla.$ Since
$\delta$ is identically zero on $\nabla$ it follows that
$\delta\equiv 0$ on $L^{\infty}(\Omega).$ Since $\delta$ commutes
with the mixing operation and every element $\lambda\in L^{0}(\Omega)$ can
be represented as
$\lambda=\sum\limits_{\alpha}\pi_{\alpha}\lambda_{\alpha},$ where
$\{\lambda_\alpha\}\subset L^{\infty}(\Omega),$ and
$\{\pi_{\alpha}\}$ is a partition of unit in  $\nabla,$ we have
$\delta(\lambda)=\delta(\sum\limits_{\alpha}\pi_{\alpha}\lambda_{\alpha})=
\sum\limits_{\alpha}\pi_{\alpha}\delta(\lambda_{\alpha})=0,$ i.e.
$\delta\equiv 0$ on $L^{0}(\Omega).$ This contradiction shows that the set
$\{\delta(\lambda): \lambda\in L^{0}(\Omega), |\lambda|\leq\textbf{1}\}$
is not order bounded in  $L^{0}(\Omega).$ Further, since  $\delta$
commutes with the mixing operations and the set  $\{\lambda:
\lambda\in L^{0}, |\lambda|\leq\textbf{1}\}$ is cyclic, the set
$\{\delta(\lambda): \lambda\in L^{0}(\Omega), |\lambda|\leq\textbf{1}\}$
is also cyclic.
 By \cite[Proposition 3]{Kud} there exist a sequence
  $\{\lambda_n\}_{n=1}^{\infty}$ in $L^{\infty}(\Omega)$ with $|\lambda_n|\leq\textbf{1}$ and an idempotent
  $\pi\in \nabla,\,\pi\neq 0,$ such that $|\delta(\lambda_n)|\geq n\pi,\,n\in \mathbb{N}.$
  The proof is complete. $\blacksquare$

Now we are in position to consider derivations on the algebra of locally measurable
operators for type I$_{\infty}$ von Neumann algebras.

\textbf{Theorem 2.6.} \emph{If $M$ is a type $I_{\infty}$ von Neumann algebra, then any
derivation on the algebra $LS(M)$ is inner.}

Proof. Since  $M$ is of type I$_{\infty}$  there exists a sequence
of  mutually orthogonal and  mutually equivalent abelian projections $\{p_n\}_{n=1}^{\infty}$ in $M$ with the central cover $\textbf{1}$
 (i.e. faithful projections).

For any bounded sequence  $\Lambda=\{\lambda_k\}$ in $Z$
define an operator  $x_\Lambda$ by
$$x_\Lambda=
\sum\limits_{k=1}^{\infty}\lambda_k p_k.$$ Then
\begin{equation}
x_\Lambda p_n=p_n x_\Lambda=\lambda_n p_n.
\end{equation}

Let  $D$ be a derivation on $LS(M),$ and let $\delta$
be its restriction onto the center of  $LS(M),$
identified with $L^{0}(\Omega).$

Take any $\lambda\in L^{0}(\Omega)$ and $n\in \mathbb{N}.$ From the identity
$$D(\lambda p_n)=D(\lambda)p_n+\lambda D(p_n)$$
 multiplying it by $p_n$ from  both sides we obtain
$$p_nD(\lambda p_n)p_n=p_n D(\lambda)p_n+\lambda p_n D(p_n)p_n.$$
Since  $p_n$ is a projection, one has that  $p_n D(p_n)p_n=0,$ and
since $D(\lambda)=\delta(\lambda)\in L^{0}(\Omega),$ we have
\begin{equation}
p_nD(\lambda p_n)p_n=\delta(\lambda)p_n.
\end{equation}

Now from the identity
$$D(x_\Lambda p_n)=D(x_\Lambda)p_n+x_\Lambda D(p_n),$$
in view of (3)  one has similarly
$$p_nD(\lambda_n p_n)p_n=p_n
D(x_\Lambda)p_n+\lambda p_n D(p_n)p_n,$$ i.e.
\begin{equation}
p_n D(\lambda_n p_n)p_n=p_nD(x_\Lambda)p_n.
\end{equation}
 (4) and (5) imply
$$p_nD(x_\Lambda)p_n=\delta(\lambda_n)p_n.$$
Further for the center-valued norm $\|\cdot\|$
on $LS(M)$ (see Section 1) we have :
 $$\|p_nD(x_\Lambda)p_n\|\leq\|p_n\|\|D(x_\Lambda)\|\|p_n\|=
\|D(x_\Lambda)\|$$ and
$$\|\delta(\lambda_n)p_n\|=|\delta(\lambda_n)|.$$
Therefore
$$\|D(x_\Lambda)\|\geq|\delta(\lambda_n)|$$
for any bounded sequence $\Lambda=\{\lambda_n\}$ in $Z.$

If we suppose that $\delta\neq 0$ then  by Lemma 2.5 there exist a
bounded sequence $\Lambda=\{\lambda_n\}$ in $Z$ and an
idempotent $\pi\in \nabla,\,\pi\neq0,$ such that
$$|\delta(\lambda_n)|\geq n\pi $$
for any $n\in \mathbb{N}.$ Thus
\begin{equation}
\|D(x_\Lambda)\|\geq n\pi
\end{equation}
for all $n\in \mathbb{N},$ i.e. $\pi=0$ -- that is  a
contradiction. Therefore $\delta\equiv 0,$ i.e. $D$ is identically
zero on the center of $LS(M),$ and therefore it is $Z$-linear. By Theorem 2.1 $D$ is inner. The proof is complete. $\blacksquare$

We shall now consider derivations on the
algebra $LS(M)$ of locally measurable operators with respect to
an arbitrary  type I von Neumann algebra $M.$

Let $M$ be a type  I von Neumann algebra. There exists a central
projection  $z_0\in M$ such that

a) $z_0M$ is a finite von Neumann algebra;

b) $z_0^{\bot}M$ is a von Neumann algebra of type  I$_{\infty}.$

Consider a derivation  $D$ on  $LS(M)$ and let $\delta$
be its restriction onto its center  $Z(S).$ By Theorem 2.6 $z_0^{\bot}D$ is inner and thus  we have $z_0^{\bot}\delta\equiv 0,$ i.e. $\delta=z_0\delta.$

Let   $D_\delta$ be the derivation on  $z_0LS(M)$ defined
as in (2) and consider its extension $D_\delta$ on $LS(M)=z_0LS(M)\oplus z_0^{\bot}LS(M)$ which is
defined as
\begin{equation}
D_\delta(x_1+x_2):=D_\delta(x_1),\,x_1\in z_0LS(M),x_2\in z_0^{\bot}LS(M).
\end{equation}

The following theorem is the main result of this section, and gives
the general form of derivations on the algebra $LS(M).$

\textbf{Theorem  2.7.} \emph{Let  $M$ be a type  $I$ von Neumann
algebra. Each
derivation  $D$ on  $LS(M)$ can be uniquely represented in
the form
$$D=D_{a}+D_{\delta}$$
where  $D_{a}$ is an inner derivation implemented by an element
$a\in LS(M),$ and $D_{\delta} $ is a derivation of the form
(7), generated by a derivation $\delta$ on the center of $LS(M)$.}

Proof. It immediately follows from Lemma  2.3 and Theorem  2.6. $\blacksquare$

\begin{center} \textbf{3. Derivations on the algebra $S(M)$}\end{center}

In this section we describe derivations on the algebra $S(M)$
of measurable operators affiliated with a type I von Neumann algebra $M.$

Let  $M$ be a type  I von Neumann algebra and let
$\mathcal{A}$ be an arbitrary subalgebra of $LS(M)$ containing
$M.$ Consider a derivation  $D:\mathcal{A}\rightarrow LS(M)$
and let us show that  $D$ can be extended to a derivation
$\tilde{D}$ on the whole $LS(M).$

Since $M$ is a type I, for an arbitrary element  $x\in LS(M)$ there exists a sequence  $\{z_n\}$
of mutually orthogonal central projections with
$\bigvee\limits_{n\in\mathbb{N}}z_n =\textbf{1}$ and  $z_n x\in M$
for all $n\in\mathbb{N}.$ Set
\begin{equation}
\tilde{D}(x)=\sum\limits_{n\geq1} z_n D(z_n x).
\end{equation}
Since every derivation  $D:\mathcal{A}\rightarrow LS(M)$ is
identically zero on central projections of  $M,$ the equality (8)
gives a well-defined derivation
$\tilde{D}:LS(M)\rightarrow LS(M)$ which coincides with $D$ on
$\mathcal{A}.$

In particular, if $D$ is $Z$-linear on $\mathcal{A},$ then $\tilde{D}$ is also $Z$-linear and by
Theorem 2.1 the derivation  $\tilde{D}$ is
inner on $LS(M)$ and therefore  $D$ is a spatial derivation on $\mathcal{A},$
i.~e. there exists an element $a\in LS(M)$ such that
$$D(x)=ax-xa$$
for all $x\in \mathcal{A}.$

Therefore we obtain the following

\textbf{Theorem 3.1.} \emph{Let     $M$ be a type  I von Neumann
algebra with the center $Z,$ and let   $\mathcal{A}$ be an
arbitrary subalgebra in $LS(M)$ containing  $M.$ Then any
$Z$-linear derivation  $D:\mathcal{A}\rightarrow LS(M)$  is spatial and
implemented by an element of  $LS(M).$}

\textbf{Corollary 3.2}. \emph{Let $M$ be a type I von Neumann algebra with the center $Z$ and let $D$
be a $Z$-linear derivation on $S(M)$ or $S(M, \tau).$  Then $D$ is spatial and implemented by an element of}
 $LS(M).$

We are now in position to improve the last result by showing that in fact such
derivations on $S(M)$ and $S(M, \tau)$ are inner.

Let us start by the consideration of the type I$_{\infty}$ case.

Let $M$ be a type I$_{\infty}$ von Neumann algebra with the center $Z$ identified with the algebra  $L^{\infty}(\Omega)$  and
let $\nabla$ be the Boolean algebra of projection from $L^{\infty}(\Omega).$

Denote by  $St(\nabla)$ the set of all elements $\lambda\in L^{\infty}(\Omega)$ of the form  $\lambda=\sum\limits_{\alpha}\pi_\alpha
t_\alpha,$  where $\{\pi_\alpha\}$ is a partition of the unit in  $\nabla,$ and
$\{t_\alpha\}\subset\mathbb{R}$ (so called step-functions).

 Suppose that $a\in LS(M),$ $a=a^{\ast}$ and consider the spectral family
$\{e_{\lambda}\}_{\lambda\in\mathbb{R}}$ of the operator  $a.$ For $\lambda\in St(\nabla),$
$\lambda=\sum\limits_{\alpha}\pi_\alpha t_\alpha$ put
$e_{\lambda}=\sum\limits_{\alpha}\pi_\alpha e_{t_\alpha}.$

Denote by $P_{\infty}(M)$ the family of all faithful projections $p$ from
$M$ such that $pMp$ is of type  I$_{\infty}.$

Set $$
\Lambda_{-}=\{\lambda\in St(\nabla): e_{\lambda}\in
P_{\infty}(M)\}
$$
and  $$
\Lambda_{+}=\{\lambda\in St(\nabla): e_{\lambda}^{\perp}\in
P_{\infty}(M)\}.
$$

\textbf{Lemma 3.3.} \emph{a) $\Lambda_{-}\neq \emptyset$ and $\Lambda_{+}\neq \emptyset;$}

\emph{b) the set $\Lambda_{+}$ (resp. $\Lambda_{-}$) is bounded from above (resp. from below);}

\emph{c) if $\lambda_{+}=\sup\Lambda_{+}$ (resp. $\lambda_{-}=\inf\Lambda_{-}$) then $\lambda\in\Lambda_{+}$ (resp. $\lambda\in\Lambda_{-}$) for all
$\lambda\in St(\nabla)$ with $\lambda+\varepsilon\textbf{1}\leq\lambda_{+}$ (resp. $\lambda-\varepsilon\textbf{1}\geq\lambda_{-}$) for some $\varepsilon>0.$}

\emph{d) if $\lambda_{+}\in L^{\infty}(\Omega),$ then $a\in S(M).$}

Proof. a) Take a sequence of projections $\{z_n\}$ from $\nabla$ such that $z_n a \in M$ for all $n\in \mathbb{N}.$
Then for $t_n<-\parallel z_na\parallel_{M}$ we have $e_{t_{n}}=0.$
Therefore for $\lambda=\sum z_n t_n$
one has $e_{\lambda}^{\bot}=\textbf{1},$ i.e. $\lambda\in\Lambda_{+}$ and hence
$\Lambda_{+}\neq \emptyset.$ Similarly  $\Lambda_{-}\neq \emptyset.$

b) Suppose that the element $\lambda=\sum\pi_\alpha\lambda_\alpha\in St(\nabla),$   satisfies the condition
$\pi_0\lambda\geq\pi_0||a||+\varepsilon\pi_0$ for an appropriate non zero $\pi_0\in\nabla,$ where $||\cdot||$ is the
center-valued norm on $LS(M).$ Without loss of
generality we may assume that $\pi_0=\pi_\alpha$ for some $\alpha,$ i.e.
$\pi_\alpha t_\alpha\geq\pi_\alpha||a||+\varepsilon\pi_\alpha.$ Then $t_\alpha\geq||\pi_\alpha a||_M+\varepsilon$
and therefore $\pi_\alpha e_{t_{\alpha}}=\pi_\alpha\textbf{1}.$ Thus $\pi_\alpha e_{t_{\alpha}}^{\bot}=0,$ i.e. $\lambda\notin\Lambda_+.$
Therefore  $\Lambda_+$ is bounded from above by the element $||a||.$
Similarly the set  $\Lambda_-$ is bounded from below by the element $-||a||.$

c) Put $$\lambda_{+}=\sup\Lambda_{+}$$ and $$\lambda_{-}=\inf\Lambda_{-}.$$
Take an element $\lambda\in St(\nabla)$ such that $\lambda+\varepsilon\textbf{1}\leq\lambda_+,$ where $\varepsilon>0.$
Suppose that $e_{\lambda}^{\bot}\notin P_\infty(M).$ Then $\pi_0 e_{\lambda}^{\bot}Me_{\lambda}^{\bot}$ is a finite von
Neumann algebra for some non zero $\pi_0\in\nabla.$ Without loss of generality we may assume
that $\pi_0=\pi_\alpha$ for some $\alpha,$ i.e. $\pi_\alpha e_{t_{\alpha}}^{\bot}$ is a finite projection. Then
$\pi_\alpha e_{t}^{\bot}$
is finite for all $t>t_\alpha.$
This means that $\pi_\alpha\lambda_+\leq\pi_\alpha t_\alpha.$

On the other hand multiplying by $\pi_\alpha$ the unequality $\lambda+\varepsilon\textbf{1}\leq\lambda_+$
we obtain that $\pi_\alpha t_\alpha+\pi_\alpha\varepsilon\leq\pi_\alpha\lambda_+.$
 Therefore $\pi_\alpha\varepsilon\leq 0.$ This contradiction implies that $\lambda\in \Lambda_+$ for all $\lambda\in St(\nabla)$ with
$\lambda+\varepsilon\textbf{1}\leq\lambda_+.$

d) Let $\lambda_{+}\in L^{\infty}(\Omega).$ Take a number $n\in\mathbb{N}$ such that $\lambda_+\leq n\textbf{1}.$ Then by the definition of $\lambda_+$ it follows that $e_{n+1}^{\bot}$
is a finite projection, i.e. $a\in S(M).$ The proof is complete. $\blacksquare$

 \textbf{Lemma 3.4.} \emph{ If $M$ is a type  I$_{\infty}$ von Neumann algebra then every derivation
 $D:M\rightarrow
S(M)$ has the form
$$
D(x)=ax-xa,\quad x\in M
$$ for an appropriate}
 $a\in S(M).$

Proof. By the Remark 1 $D$ maps the center $Z$ of $M$ into the center of  $S(M)$  which coincides with $Z$
by Proposition 1.2, i.e. we obtain a derivation $D$
on commutative von Neumann algebra $Z.$ Therefore $D|_{Z}=0.$ Thus $D(\lambda x)=D(\lambda)x+\lambda D(x)=\lambda D(x)$ for all
$\lambda\in Z,$ i.e. $D$ is $Z$-linear.

By Theorem 3.1 there exists an element $a\in LS(M)$  such that   $D(x)=ax-xa$ for all
$x\in M.$

Let us prove that one can choose the element $a$ from $S(M).$

For  $x\in M$ we have
 $$
 (a+a^{\ast})x-x(a+a^{\ast})=(ax-xa)-(ax^{\ast}-x^{\ast}a)^{\ast}=D(x)-D(x^{\ast})^{\ast}\in S(M)
 $$
and
 $$
 (a-a^{\ast})x-x(a-a^{\ast})=D(x)+D(x^{\ast})^{\ast}\in S(M).
 $$
This means that the elements  $a+a^{\ast}$ and $a-a^{\ast}$ implement derivations from $M$
into  $S(M).$ Since $a=\frac{\textstyle
a+a^{\ast}}{\textstyle 2}+i\frac{\textstyle a-a^{\ast}}{\textstyle
2i},$ it is sufficient to consider the case where  $a$ is a self-adjoint element.

Consider the elements  $\lambda_{+},\,\lambda_{-}\in L^{0}$ defined in Lemma 3.3 c) and let us prove that
$\lambda_{+}-\lambda_{-}\in L^{\infty}(\Omega).$
Lemma 3.3 c) implies that there exists an element $\lambda_1\in \Lambda_-$
such that   $\lambda_-+\frac{\textstyle
1}{\textstyle 16}\leq\lambda_1\leq\lambda_--\frac{\textstyle 1}{\textstyle
2}.$ Since $D(x)=(a-\lambda_1)x-x(a-\lambda_1),$ replacing $a$ by $a-\lambda_1,$
 we may assume that $\lambda_-\leq\frac{\textstyle 1}{\textstyle
8}.$ Then  $e_\varepsilon\in P_\infty(M)$ for all $\varepsilon>\frac{\textstyle 1}{\textstyle
8}.$

Suppose that $\lambda_+\notin L^{\infty}(\Omega).$ Passing if necessary to the subalgebra $zM,$ where $z$ is a non zero central
projection in $M$ with $z\lambda_+\geq z,$ we may assume without loss of generality that $\lambda_+\geq\textbf{1}.$

First let us consider the particular case where $M$ is of type I$_{\aleph_{0}},$ where $\aleph_{0}$
 is the countable cardinal number. Take an element $\lambda_0\in St(\nabla)$
such that $\lambda_+-\frac{\textstyle 1}{\textstyle 2}\leq\lambda_0\leq\lambda_+-\frac{\textstyle 1}{\textstyle
4}.$ By Lemma 3.3 c)
we have $e_{\lambda_{0}}^{\bot}\in P_\infty(M).$
Since $e_{\lambda_{0}}^{\bot}Me_{\lambda_{0}}^{\bot}$ and $e_{\frac{1}{4}}Me_{\frac{1}{4}}$
are algebras of type I$_{\aleph_{0}},$ the projections $p_1=e_{\lambda_{0}}^{\bot}$ and $p_2=e_{\frac{1}{4}}$ are equivalent.
From $\lambda_0e_{\lambda_{0}}^{\bot}\leq ae_{\lambda_{0}}^{\bot}$ it follows that
$\lambda_0p_1\leq p_1ap_1.$
 Since $p_1Mp_1$ is of type I$_{\aleph_{0}},$
the center of the algebra $S(p_1Mp_1)$ coincides with the center of the algebra $p_1Mp_1$ (Proposition 1.2)
and therefore $\lambda_0p_1\notin S(p_1Mp_1),$ because $\lambda_0p_1$ is a central unbounded element in $LS(p_1Mp_1).$
 Therefore $ap_1=p_1ap_1\notin S(p_1Mp_1).$

 Let  $u$ be a partial isometry in  $M$ such that
 $uu^{\ast}=p_1,$ $u^{\ast}u=p_2.$ Put  $p=p_1+p_2.$
Consider the derivation  $D_1$ from  $pMp$ into $pS(M)p=
S(pMp)$ defined as
$$
D_1(x)=pD(x)p,\, x\in  pMp.
$$
This  derivation is implemented by the element $ap=pap,$ i.e.
 $$
 D_1(x)=apx-xap,\, x\in pMp.
 $$
Since $ap_2\in pMp,$ the element $b=ap_1=ap-ap_2$
 implements a derivation  $D_2$ from $pMp$ into $S(pMp).$

Since  $D_2(u+u^{\ast})=b(u+u^{\ast})-(u+u^{\ast})b,$
it follows that $b(u+u^{\ast})-(u+u^{\ast})b\in S(M).$
From  $up_1=p_1 u^{\ast}=0$ it follows that
  $bu-u^{\ast}b\in S(M).$ Multiplying this by  $u$ from the left side we obtain
 $ubu-uu^{\ast}b\in S(M).$ From $ub=0,\,uu^{\ast}=p_1,$ it follows that $p_1 b\in S(M),$ i.e.
 $ap_1\in S(M).$
This contradicts the above relation   $ap_1\notin S(M).$ The contradiction shows that  $\lambda_{+}\in L^{\infty}(\Omega).$
Now Lemma 3.3 d) implies that $a\in S(M).$

Let us consider the case of general type I$_\infty$ von Neumann algebra $M.$ Take an element
$\lambda_0\in St(\nabla)$
such that $\lambda_+-\frac{\textstyle 1}{\textstyle 2}\leq\lambda_0\leq\lambda_+-\frac{\textstyle 1}{\textstyle
4}.$ Lemma 3.3 c)  implies that  $e_{\lambda_{0}}^{\bot}\in P_\infty(M).$
Consider projections $p_1$ and $p_2$ with the central cover $\textbf{1}$
 such that $p_1\leq e_{\lambda_{0}}^{\bot},\,p_2\leq e_{\frac{1}{4}}$
and such that $p_iMp_i$ are of type I$_{\aleph_{0}},$ $i=1,2.$
 Put $p=p_1+p_2.$ Consider the derivation $D_p$ from $pMp$ into
$pS(M)p$ defined as
$$
D_p(x)=pD(x)p,\, x\in  pMp.
$$
Since $pMp$ is of type I$_{\aleph_{0}}$ the above case implies that $pap\in S(M)$ and therefore $p_1ap_1\in S(M).$
On the other hand   $\lambda_0p_1\leq p_1ap_1$ and $\lambda_0 p_1\notin S(M).$ From this contradiction it follows that
 $\lambda_{+}\in L^{\infty}(\Omega).$
By Lemma 3.3 d) we obtain that $a\in S(M).$ The proof is complete. $\blacksquare$

From the above results we obtain

\textbf{Lemma 3.5.} \emph{Let $M$ be a type I von Neumann algebra with the center $Z.$ Then
every $Z$-linear derivation $D$ on the algebra $S(M)$  is inner. In particular, if $M$ is a type I$_\infty$
then every derivation on $S(M)$ is inner.}

Now let $M$ be an arbitrary  type  I von Neumann algebra and let $z_0$ be  the central
projection  in $M$ such that $z_0M$ is a finite von Neumann algebra and
 $z_0^{\bot}M$ is a von Neumann algebra of type  I$_{\infty}.$
Consider a derivation  $D$ on  $S(M)$ and let $\delta$
be its restriction onto its center  $Z(S).$ By Lemma 3.5 the derivation $z_0^{\bot}D$ is inner and thus  we have $z_0^{\bot}\delta\equiv 0,$ i.e. $\delta=z_0\delta.$

Since $z_0M$ is a finite type I  von Neumann algebra, we have that $z_0LS(M)=z_0S(M).$ Let   $D_\delta$ be the derivation on  $z_0S(M)=z_0LS(M)$ defined
as in (2).

Finally Lemma 2.3 and Lemma 3.5 imply the following main result the present section.

\textbf{Theorem 3.6.} \emph{Let $M$ be a type  I von Neumann algebra. Then every dereivation $D$
on the algebra  $S(M)$  can be uniquely represented in the form}
$$D=D_a+D_\delta,$$
\emph{where  $D_a$ is inner and implemented by an element $a\in S(M)$
and  $D_\delta$ is the derivation of the form (7) generated by a derivation $\delta$ on the center of} $S(M).$

\begin{center} \textbf{4. Derivations on the algebra $S(M,
\tau)$}\end{center}

   In this section we present a general
form of derivations on the algebra  $S(M,\tau)$ of $\tau$-measurable
operators affiliated with a type I von Neumann algebra $M$ and a faithful normal semi-finite trace $\tau.$

\textbf{Theorem 4.1.} \emph{Let $M$ be a type I von Neumann algebra with the center $Z$ and a faithful normal
semi-finite trace $\tau.$ Then
every $Z$-linear derivation $D$ on the algebra $S(M, \tau)$ is inner.
  In particular, if $M$ is a type I$_\infty$
then every derivation on $S(M, \tau)$ is inner.}

Proof. By Theorem 3.1 $D(x)=ax-xa$ for some $a\in
LS(M)$ and all $x\in S(M, \tau).$  Let us show that the element   $a$  can be chosen from the algebra $S(M, \tau).$

\emph{Case 1.}  $M$ is a homogeneous type
I$_n,$ $n\in\mathbb{N}$ von Neumann algebra. Then
$LS(M)=S(M)\cong M_{n}(L^{0}(\Omega)).$ As in Lemma 3.3 we may assume that
 $a=a^{\ast}.$ By  \cite[Theorem 3.5]{Kus3} *-isomorphism between $S(M)$ and $M_{n}(L^{0}(\Omega))$ can be a chosen
 such that
the element  $a$ can be represented as
$a=\sum\limits_{i=1}^{n}\lambda_{i}e_{i,i},$  where $\lambda_{i}=\overline{\lambda_{i}}\in
L^{0}(\Omega),$ $i=\overline{1, n},$
 $\lambda_{1}\geq\cdots\geq\lambda_{n}.$

 Put  $u=\sum\limits_{j=1}^{n}e_{j,n-j+1}.$
Then   $$
  D_{a}(u)=au-ua=\sum\limits_{i=1}^{n}(\lambda_{i}-\lambda_{n-i+1})e_{i,n-i+1}
  $$
and
 $$
  D_{a}(u)^{\ast}=\sum\limits_{i=1}^{n}(\lambda_{i}-\lambda_{n-i+1})e_{n-i+1,i}.
  $$
Therefore
$D_{a}(u)^{\ast}D_{a}(u)=\sum\limits_{i=1}^{n}(\lambda_{i}-\lambda_{n-i+1})^{2}e_{i,i},$
and thus
$|D_{a}(u)|=\sum\limits_{i=1}^{n}|\lambda_{i}-\lambda_{n-i+1}|e_{i,i},$
Since $\lambda_{1}\geq\cdots\geq\lambda_{n},$ we have
\begin{equation}
|\lambda_{i}-\lambda_{n-i+1}|\geq |\lambda_{i}-\lambda_{[\frac{
n+1}{2}]}|
\end{equation}
  for all   $i\in\overline{1, n}.$

Set $b=\sum\limits_{i=1}^{n}|\lambda_{i}-\lambda_{[\frac{n+1}{2}]}|e_{i,i}.$
From (9) we obtain that  $|D_{a}(u)|\geq b,$ and thus  $b\in S(M, \tau).$

 Set  $v=\sum\limits_{i=1}^{[\frac{n+1}{2}]}e_{i, i}-\sum\limits_{j=[\frac{n+1}{2}]}^{n}e_{j,j}.$ Then
$vb=a-\lambda_{[\frac{n+1}{2}]}\textbf{1}$ and  $vb\in S(M, \tau).$
Therefore $a-\lambda_{[\frac{n+1}{2}]}\textbf{1}\in S(M, \tau)$ and this element also implements the derivation $D_a.$

 \emph{Case 2.} Let $M$ be a finite type  I von Neumann algebra. Then
$$
LS(M)=S(M)\cong\prod_{n\in F}M_n(L^{0}(\Omega_n),
$$
where  $F\subseteq\mathbb{N}.$   Therefore   $a=\{a_{n}\},$ where
$a_n=\sum\limits_{i=1}^{n}\lambda_{i}^{(n)}e_{i,i}^{(n)},$
$\lambda_{1}^{(n)}\geq\cdots\geq\lambda_{n}^{(n)},$
$\lambda_{i}^{(n)}\in L^{0}(\Omega_n)$ and  $e_{i,j}^{(n)}$ are the matrix units in  $M_{n}(L^{0}(\Omega_n)),$ $i, j=\overline{1,
n},$ $n\in F.$

For each $n\in F$ consider the following elements in $M_n(L^{0}(\Omega_n)$
$$
b_n=\sum\limits_{i=1}^{n}|\lambda_{i}^{(n)}-\lambda_{[\frac{n+1}{2}]}^{(n)}|e_{i,i}^{(n)}
$$
and
$$v_n=\sum\limits_{i=1}^{[\frac{n+1}{2}]}e_{i, i}^{(n)}-\sum\limits_{j=[\frac{n+1}{2}]}^{n}e_{j,j}^{(n)}.
$$
Set   $b=\{b_{n}\}_{n\in F}$ and $v=\{v_{n}\}_{n\in F}.$ Consider the element
    $$\lambda=\{\lambda_{[\frac{n+1}{2}]}\}_{n\in F}\in L^{0}(\Omega)\cong\prod\limits_{n\in F} L^{0}(\Omega_n).$$
Similar to the case 1 we obtain that $a-\lambda\textbf{1}=vb\in S(M, \tau).$

\emph{Case  3.}  $M$ is a type  I$_{\infty}$ von Neumann algebra. Since
 $S(M, \tau)\subseteq S(M)$ by Lemma 3.4 there exists an element $a\in S(M)$
such that $D(x)=ax-xa$ for all $x\in M.$ Let us show that $a$ can be picked from the algebra $S(M, \tau).$
Since $a\in S(M),$ there exists $\lambda\in\mathbb{R},\,\lambda>0$ such that $e_{\lambda}^{\perp}$ is a finite projection. Then
$e_{\lambda}Me_{\lambda}$ is a type  I$_{\infty}$ von Neumann algebra and thus there exists a projection
 $q\leq e_{\lambda}$ such that  $q\sim p.$ Let
$u$ be a partial isometry in $M$ such that $uu^{\ast}=p,$
$u^{\ast}u=q.$ Similar to Lemma 3.4 we obtain that  $uapu-uu^{\ast}ap\in
S(M, \tau)$ and   $ap\in S(M, \tau).$ Therefore   $a\in S(M,
\tau).$ The proof is complete. $\blacksquare$

Let $N$ be a commutative von Neumann algebra, then $N\cong
L^{\infty}(\Omega)$ for an appropriate measure space
$(\Omega,\Sigma,\mu).$ It has been proved in  \cite{Ber},
\cite{Kus1} that the algebra $LS(N)=S(N)\cong L^{0}(\Omega)$
admits non trivial derivations if and only if the measure space
$(\Omega,\Sigma,\mu)$ is not atomic.

Let  $\tau$ be a faithful normal semi-finite trace on the commutative von
Neumann algebra $N$ and suppose that the Boolean
algebra $P(N)$ of projections is not atomic. This means that there exists a projection  $z\in N$ with
$\tau(z)<\infty$ such that the Boolean algebra of projection in $zN$ is continuous (i.e. has no atom). Since
$zS(N, \tau)=zS_0(N, \tau)=zS(N)=S(zN),$ the algebra $zS(N, \tau)$ (resp. $zS_0(N, \tau)$) admits a non trivial derivation  $\delta.$ Putting
$$\delta_0(x)=\delta(zx),\, x\in S(N, \tau)$$
we obtain a non trivial derivation $\delta_0$ on the algebra $S(N, \tau).$ Therefore,
we have that if a commutative von Neumann
algebra $N$
has a non atomic Boolean algebra of projections then the algebra  $S(N, \tau)$ admits a non zero derivation.

Given an arbitrary derivation $\delta$ on $S(N,\tau)$ or $S_{0}(N,\tau)$ the element
$$z_\delta=\inf\{z\in P(N): z\delta=\delta\}$$
is called the support of the derivation  $\delta.$

\textbf{Lemma 4.2.} \emph{If $N$ is a commutative von Neumann algebra with a
faithful normal semi-finite trace  $\tau$
and $\delta$  is a derivation on $S(N,\tau)$ or $S_{0}(N,\tau)$, then $\tau(z_{\delta})<\infty.$}

Proof. Let us give proof for the algebra $S_0(N, \tau),$ since the case of $S(N, \tau)$ is similar and simpler.
Suppose the opposite, i.e.  $\tau(z_{\delta})=\infty.$ Then
there exists a sequence of mutually orthogonal projections
$z_{n}\in N,\,n=1,2...,$ with $z_{n}\leq z_{\delta},\ 1\leq
\tau(z_{n})<\infty.$ For $z=\sup\limits_{n}{z_{n}}$ we have
$\tau(z)=\infty.$ Since $\tau(z_{n})<\infty$  for all $n=1,2...,$
it follows that $z_{n}S_{0}(N,\tau)=z_n S(N)=S(z_{n}N).$
Define a derivation $\delta_{n}: S(z_{n}N)\rightarrow S(z_{n}N)$ by
$$\delta_{n}(x)=z_{n}\delta(x), \ x\in S(z_{n}N).$$
Since  $z_{\delta_{n}}=z_{n}$,
Lemma 3.5  implies that for each $n\in\mathbb{N}$
there exists an element $\lambda_{n}\in z_{n}N$ such that
$|\lambda_{n}|\leq n^{-1}z_{n}$ and $|\delta_{n}(\lambda_{n})|\geq
z_{n}.$

Put $\lambda=\sum\limits_{n\geq 1}\lambda_{n}. $
Then $|\lambda|\leq\sum\limits_{n\geq 1}n^{-1}z_n$ and therefore $\lambda\in S_{0}(N,
\tau).$ On other hand
$$|\delta(\lambda)|=|\delta(\sum\limits_{n\geq 1}\lambda_{n})|=|\delta(\sum\limits_{n\geq 1}z_n\lambda_{n})|=
|\sum\limits_{n\geq 1}z_n \delta(\lambda_{n})|=\sum\limits_{n\geq 1}|\delta_n(\lambda_{n})|\geq \sum\limits_{n\geq 1}z_n=z,$$
i.e. $|\delta(\lambda)|\geq z.$ But $\tau(z)=\infty$, i.e. $z\notin S_{0}(N, \tau).$ Therefore
$\delta(\lambda)\notin S_{0}(N, \tau).$
The contradiction shows that $\tau(z_{\delta})<\infty.$ The proof is complete. $\blacksquare$

Let  $M$ be a homogeneous von Neumann algebra of type I$_{n},
n \in \mathbb{N}$,  with the center $Z$ and a faithful normal
semi-finite trace  $\tau.$ Then  the algebra $M$ is *-isomorphic
with the algebra $M_n(Z)$ of all  $n\times n$- matrices over $Z,$
and the algebra  $S(M, \tau)$ is *-isomorphic with the algebra
 $M_n(S(Z, \tau_Z))$ of all  $n\times n$ matrices over
$S(Z, \tau_Z),$ where  $\tau_Z$ is the restriction of the trace
 $\tau$ onto the center  $Z.$

Now let  $M$ be an arbitrary finite von Neumann algebra of type I
with the center $Z$ and  let $\{z_n\}_{n\in F},$
$F\subseteq\mathbb{N},$ be a  family of central projections from $M$ with
$\sup\limits_{n\in F}z_n=\textbf{1}$ such that the algebra  $M$ is
*-isomorphic with the  $C^{*}$-product of von Neumann algebras
$z_n M$ of type  I$_{n}$ respectively, $n\in F,$ i.e.
$$M\cong\bigoplus\limits_{n\in F}z_n M.$$
In this  case we have that
$$S(M, \tau)\subseteq\prod\limits_{n\in F}S(z_n M, \tau_n),$$
where  $\tau_n$ is the restriction of the trace $\tau$ onto $z_n
M,\,n\in F.$

Suppose that    $D$ is a derivation on  $S(M, \tau),$ and let
$\delta$ be its restriction onto the center  $S(Z, \tau_Z).$
Since  $\delta$ maps each  $z_nS(Z, \tau_Z)\cong Z(S(z_n M,
\tau_n))$ into itself,  $\delta$ generates a derivation $\delta_n$
on $z_nS(Z, \tau_Z)$ for each $n\in F.$

Let     $D_{\delta_n}$ be the derivation on the matrix algebra
$M_n(z_nZ(S(M, \tau)))\cong S(z_nM, \tau_n)$ defined as in
(1). Put
\begin{equation}
D_\delta(\{x_n\}_{n\in F})=\{D_{\delta_n}(x_n)\},\,\{x_n\}_{n\in F}\in S(M, \tau).
\end{equation}
By Lemma 4.2 $\tau(z_\delta)<\infty,$ thus
$$z_\delta S(M, \tau)=z_\delta S(M)\cong z_\delta\prod\limits_{n\in F}S(z_n M)=
z_\delta\prod\limits_{n\in F}S(z_n M, \tau_n),$$
 and therefore $\{D_{\delta_n}(x_n)\}\in z_\delta S(M, \tau)$ for all
  $\{x_n\}_{n\in F}\in S(M, \tau).$ Hence we obtain that  the map  $D$  is a derivation on $S(M, \tau).$

Similar to Lemma  2.3 one can prove the following.

\textbf{Lemma 4.3.} \emph{Let  $M$ be a finite von Neumann algebra
of type I with a faithful normal semi-finite trace $\tau.$ Each
derivation  $D$ on the algebra  $S(M, \tau)$ can be uniquely
represented in the form
$$D=D_{a}+D_{\delta ,}$$
where  $D_{a}$ is an inner derivation implemented by an element
$a\in S(M, \tau),$ and $D_{\delta} $ is a derivation given as
(10).}

Finally Theorem 4.1 and Lemma 4.3 imply the following main result the present section.

\textbf{Theorem 4.4.} \emph{Let $M$ be a type  I von Neumann algebra with a faithful normal semi-finite trace $\tau.$ Then every derivation $D$
on the algebra    $S(M, \tau)$ can be uniquely represented in the form}
$$D=D_a+D_\delta,$$
\emph{where  $D_a$ is inner and implemented by an element  $a\in S(M, \tau)$
and  $D_\delta$ is the derivation of the form (10) generated by a derivation $\delta$ on the center of}
 $S(M, \tau)$.

 If we consider the measure topology $t_{\tau}$ on the algebra $S(M, \tau)$ (see Section 1) then it is clear that every non-zero derivation of the form $D_\delta$ is discontinuous in $t_\tau.$ Therefore the above Theorem 4.4 implies

 \textbf{Corollary 4.5.} \emph{Let $M$ be a type  I von Neumann algebra with a faithful normal semi-finite trace $\tau.$ A derivation $D$
on the algebra    $S(M, \tau)$ is inner if and only if it is continuous in the measure topology.}

\begin{center} \textbf{5. Derivations on the algebra $S_0(M,
\tau)$ }\end{center}

In this section we  describe derivations on the algebra
$S_0(M, \tau)$ of all $\tau$-compact operators for type I von
Neumann algebra $M$ with a faithful normal semi-finite trace $\tau.$

It should be noted that the centers of the algebras $LS(M),$ $S(M)$ and $S(M, \tau)$
for general von Neumann algebra  $M$ contain $Z.$ This was an essential
 point in the proof of theorems
concerning the description of derivations on these algebras.
Proposition 1.2 shows that this is not the case
for the algebra  $S_{0}(M,\tau)$ because the center of this
algebra may be trivial. Thus the methods of previous sections can not be directly applied for the description of derivations of
the algebra $S_0(M, \tau).$

First recall the following main result of the paper  \cite{Alb2}.

\textbf{Theorem 5.1.} \emph{Let $M$ be a type I von Neumann algebra with the center $Z$ and a faithful normal
semi-finite trace $\tau.$ Then
every $Z$-linear derivation $D$ on the algebra $S_0(M, \tau)$ is spatial and implemented by an element from $S(M, \tau).$}

The main result of this section will be proved  step by step in several particular cases.

For a finite type I von Neumann algebras
we have

 \textbf{Lemma 5.2.} \emph{Let $M$ be a finite von Neumann
algebra of type I with the center  $Z$ and let
$D:S_{0}(M,\tau)\rightarrow S_{0}(M,\tau)$ be a derivation. If
$D(\lambda)=0$ for every $\lambda$ from the center
$Z(S_{0}(M,\tau))$ of $S_{0}(M,\tau),$ then $D$ is $Z$-linear.}

Proof. Take  $\lambda\in Z$ and choose a central projection $z$ in
$M$ with  $\tau(z)<\infty.$ Since  $z,\,z\lambda\in
Z(S_{0}(M,\tau)),$ we have that  $D(z)=D(z\lambda)=0.$

For  $x\in S_{0}(M,\tau)$ one has  $$D(z\lambda x)=D(z\lambda)x+z
\lambda D(x)=z\lambda D(x),$$ i.e.
$$D(z\lambda
x)=z\lambda D(x).$$ On the other hand
$$D(z\lambda x)=D(z)\lambda x+zD(\lambda x)=zD(\lambda x),$$
i.e.
$$D(z\lambda
x)=zD(\lambda x).$$ Therefore $zD(\lambda x)=z\lambda D( x). $
Since  $z$ is an arbitrary with  $\tau(z)<\infty$ this implies
(taking $z\uparrow \textbf{1}$) that $D(\lambda x)=\lambda D( x)$
for all $\lambda\in Z$ and $x\in S_0(M, \tau),$ i.e. $D$ is
$Z$-linear. The proof is complete. $\blacksquare$

 Now let $M$ be a type I$_n$  von Neumann algebra with a finite trace $\tau.$ Then  $S_0(M, \tau)=S(M, \tau)=S(M).$ Consider a family
    $\{e_i\}_{i=1}^{n}$ of mutually orthogonal and mutually equivalent abelian projections in
    the von Neumann algebra  $M.$
Put  $e=\sum\limits_{i=1}^{n-1}e_i.$ Then  $eMe$ is a von Neumann
algebra of type  I$_{n-1},$ and
$$S_0(Z, \tau_Z)\cong Z(eS_0(M, \tau)e)\cong Z(S_0(M, \tau)).$$

\textbf{Remark 2.} From now on we shall identify these isomorphic
abelian von Neumann algebras. In this case the element $\lambda$ from $S_0(Z, \tau_Z)$  corresponds to $\lambda e$
from $Z(eS_0(M, \tau)e)$ and to $\lambda\textbf{1}$ from
$Z(S_0(M, \tau)).$

Consider a derivation $D$ on the algebra  $S_0(M, \tau).$ Since
$D$ maps   $Z(S_0(M, \tau))$ into itself, its
restriction  $D|_{Z(S_0(M, \tau))}$ induces a derivation
 $\delta$ on $S_0(Z, \tau_Z)\cong  Z(S_0(M, \tau)),$ i.e.
 $$D(\lambda\textbf{1})=\delta(\lambda)\textbf{1},\,\lambda\in S_0(Z, \tau_Z).$$
 Let  $D_e$ be the derivation on $eS_0(M, \tau)e$ defined as
 $$D_e(x)=eD(x)e, \,x\in eS_0(M, \tau)e.$$
 Since  $Z(eS_0(M, \tau)e)\cong Z(S_0(M, \tau)),$ the restriction of
 $D_e$ onto   $Z(eS_0(M, \tau)e)$ also generates a derivation, denoted by
  $\delta_e,$ on $S_0(Z, \tau_Z),$ i.e.
 $$D_e(\lambda e)=\delta_e(\lambda)e,\,\lambda\in S_0(Z, \tau_Z).$$

\textbf{Lemma 5.3.} \emph{The derivations $\delta$ and $\delta_e$ on
$S_0(Z, \tau_Z)$ coincide.}

 Proof.  Since  $e$ is a projection it is clear that  $eD(e)e=0$ and therefore
  $$\delta_e(\lambda)e=D_e(\lambda e)=eD(\lambda e)e=
 eD(\lambda\textbf{1})e+e\lambda D(e)e=eD(\lambda\textbf{1})e=\delta(\lambda)e,$$ i.e.
 $$\delta_e(\lambda)e=\delta(\lambda)e$$
 for any  $\lambda\in S_0(M, \tau).$ Therefore (see Remark 2)
 $\delta_e(\lambda)=\delta(\lambda),$ i.e.
  $\delta\equiv \delta_e.$ The proof is complete. $\blacksquare$

Now similar to the proof of Lemma 2.2 from  Lemma 5.2 we obtain following results which describes derivations on the algebra of
$\tau$-compact operators for type  I$_n,$ $n\in\mathbb{N},$ von
Neumann algebras.

 \textbf{Lemma 5.4.} \emph{Let  $M$ be a homogenous von Neumann algebra of type
  $I_{n}, n \in \mathbb{N}$, with a faithful normal semi-finite
  trace $\tau.$ Every derivation  $D$ on the algebra $S_0(M, \tau)$ can be uniquely represented as a sum
  $$D=D_{a}+D_{\delta ,}$$ where  $D_{a}$ is a spatial derivation implemented by an element  $a\in S(M, \tau)$
while $D_{\delta} $ is the derivation of the form (1) generated by
a derivation $\delta$ on the center of $S_0(M, \tau)$ identified
with $S_0(Z, \tau_Z)$.}

We are now in position to prove one of the main results of this section.

 \textbf{Theorem 5.5.} \emph{If $M$ is a type  $I_\infty$ von Neumann algebra with a faithful normal semi-finite trace
  $\tau,$ then every derivation on the algebra $S_0(M, \tau)$ is spatial and implemented by an element of the algebra
   $S(M, \tau).$ }

 The proof of the theorem consists of several lemmata.

\textbf{Lemma  5.6.} \emph{Let   $z\in Z$ be a central projection
from $M$ and let  $x\in S_0(M, \tau).$ Then}
 $$D(zx)=zD(x).$$

Proof. Without loss of generality we may suppose that  $x\geq
0,$ i.e. $x=y^{2}$ for some $y\in S_0(M, \tau).$ From the Leibniz
rule for derivations we obtain
$$D(zx)=D(zyzy)=D(zy)zy+zyD(zy)=
z[D(zy)y+yD(zy)].$$ Therefore
$$z^{\perp}D(zx)=0.$$
Similarly we have that
$$zD(z^{\perp}x)=0.$$
Further
$$zD(x)=zD((z+z^{\perp})x)=zD(zx)+zD(z^{\perp}x)=zD(zx),$$
i.e.
$$zD(x)=zD(zx).$$
On the other hand
$$D(zx)=(z+z^{\perp})D(zx)=zD(zx)+z^{\perp}D(zx)=zD(zx),$$
i.e.
$$D(zx)=zD(zx).$$
 Therefore   $$D(zx)=zD(x).$$ The proof is complete. $\blacksquare$

  \textbf{Lemma  5.7.} \emph{Suppose that $\lambda\in Z,$
  $p\in P(M),\, \tau(p)<\infty.$ Put $y=D(\lambda p)-\lambda D(p).$ Then}
 $$p^{\perp}yp^{\perp}=0.$$

 Proof. From $$D(p)=D(pp)=D(p)p+pD(p)$$ and
  $$D(\lambda p)=D(\lambda pp)=D(\lambda p)p+\lambda pD(p)$$ we
  obtain
 $$p^{\perp}D(\lambda p)p^{\perp}=p^{\perp}\lambda D(p)p^{\perp}=0$$ and in particular $p^{\perp}yp^{\perp}=0.$
The proof is complete. $\blacksquare$

 \textbf{Lemma  5.8.} \emph{For each  $\lambda\in Z$ and for every abelian projection
  $p\in P(M)$ with $\tau(p)<\infty$ we have}
 $$D(\lambda p)=\lambda D(p).$$

Proof.  Let $z$ be the central cover of the projection $p.$  Lemma 5.6 implies
that the derivation $D$ maps the algebra $zS_0(M, \tau)$ into
itself. Therefore passing if necessary to the algebra $zM$ and to
the derivation  $zD$ we may assume without loss of generality
that $z=\textbf{1},$ i.e. that $p$ is a faithful projection. Take
an arbitrary faithful projection  $p_0$
such that $p_0\leq p^{\perp}$ and such that the von Neumann
algebra $p_0Mp_0$ is of type  I$_{\aleph_{0}},$ where $\aleph_{0}$
is the countable cardinal number. Then there exists a sequence of
mutually orthogonal and pairwise equivalent abelian projections
$\{p_{n}\}_{n=2}^{\infty}$ in $M$ with $\sum\limits_{n=2}^{\infty}
p_{n}=p_0.$ Putting $p_1=p$ we obtain that the projections $p_1$ and $p_n$
are equivalent ($p_1\sim p_n$) and thus
 $\tau(p_n)=\tau(p_1)<\infty$ for all  $n\geq 2.$

Set $e_{n}=\sum\limits_{k=1}^{n} p_{k},\ n\geq 1$. Then
$e_{n}Me_{n}$ is a homogeneous von Neumann algebra of type
I$_{n},$ and the restriction $\tau_{n}$ of the trace
$\tau$ onto $e_{n}Me_{n}$ is finite, and therefore  $e_{n}S_0(M,
\tau)e_{n}=S(e_{n}Me_{n}),\,n=1,2....$

Define a derivation $D_n$ on $e_{n}S_0(M, \tau)e_{n}$ as
follows
$$D_{n}(x)=e_{n}D(x)e_{n},\, x \in e_{n}S_0(M,
\tau)e_{n}.$$ By Lemma  5.4 for each  $n$ there exists an element
$a_{n} \in e_{n}S_0(M, \tau)e_{n}$ and a derivation $\delta_{n}$
on  $e_{1}S_0(M, \tau)e_{1}$ identified with $Z(e_{n}S_0(M,
\tau)e_{n})$ (see Remark 2) such that
\begin{equation}
D_{n}=D_{a_{n}}+D_{\delta_{n}}.
\end{equation}
Since  $D_n=e_n D_{n+1}e_n$ Lemma 5.3 implies that
$\delta_{n}=\delta_{n+1},\ n\geq 1.$ Denote $\delta=\delta_n.$

Given a sequence  $\Lambda=\{\lambda_{n}\}$ in  $Z$ with
$|\lambda_{n}|\leq \frac{\textstyle 1}{\textstyle n}\textbf{1},\ n
\in \mathbb{N},$ put
$$x_{\Lambda}=\sum\limits_{n=1}^{\infty}\lambda_{n}p_{n}.$$
Let us show that  $x_{\Lambda}\in S_0(M, \tau).$
For an arbitrary  $\varepsilon>0$ there exists  $n_{0}\in
\mathbb{N}$ such that $\frac{\textstyle 1}{\textstyle
n_0}<\varepsilon.$ Set
$$p_\varepsilon=\textbf{1}-\sum
\limits_{n=1}^{n_{0}-1}p_{n},$$ then
 $\tau(p^{\perp}_\varepsilon)=\tau(\sum
\limits_{n=1}^{n_{0}-1}p_{n})=(n_{0}-1)\tau(p_1)<\infty.$ Moreover
$$||x_{\Lambda}p_\varepsilon||_{M}=||\sum
\limits_{n=n_{0}}^{\infty}\lambda_{n}p_{n}||_{M}=\sup\limits_{n\geq
n_0}||\lambda_{n}||_{M}\leq {1 \over n_{0}}<\varepsilon.$$
 This means that  $x_{\Lambda}\in
S_{0}(M,\tau).$ For each  $n \in \mathbb{N}$ we have
$$x_{\Lambda}p_{n}=p_{n}x_{\Lambda}=\lambda_{n} p_{n}.$$
Similar to the proof of  (5) in Theorem 2.5 we obtain
\begin{equation}
p_{n}D(\lambda_{n}p_{n})p_{n}=p_{n}D(x_{\Lambda})p_{n}.
\end{equation}

On the other hand
$$p_{n}D(\lambda_{n}p_{n})p_{n}=p_{n}e_{n}D(\lambda_{n}p_{n})e_{n}p_{n}=p_{n}D_{n}(\lambda_{n}p_{n})p_{n}.$$
From  (11) we obtain
$$p_{n}D(\lambda_{n}p_{n})p_{n}=p_{n}D_{a_{n}}(\lambda_{n}p_{n})p_{n}+p_{n}D_{\delta}(\lambda_{n}p_{n})p_{n}.$$
Since  $D_{a_{n}}$ is a spatial derivation (and hence it is
$Z$-linear), we have that
$$p_{n}D_{a_{n}}(\lambda_{n}p_{n})p_{n}=\lambda_{n}p_{n}D_{a_{n}}(p_{n})p_{n}=0.$$
From
$$p_{n}D_{\delta}(\lambda_{n}p_{n})p_{n}=\delta(\lambda_{n})p_{n},$$
we obtain
\begin{equation}
p_{n}D(\lambda_{n}p_{n})p_{n}=\delta(\lambda_{n})p_{n}.
\end{equation}
 Now  (12) and  (13) imply
$$p_{n}D(x_{\Lambda})p_{n}=\delta(\lambda_{n})p_{n}.$$

Suppose that  $\delta\neq 0.$ Then  Lemma 2.5 implies
the existence of a sequence  $\Lambda=\{\lambda_{n}\}$ in $Z$ with
 $|\lambda_{n}|\leq \frac{\textstyle 1}{\textstyle n}\textbf{1},\
n \in \mathbb{N},$ and a projection  $\pi\in Z,\,\pi\neq0$ such
that
$$|\delta(\lambda_{n})|\geq n\pi, \ n\in \mathbb{N}.$$
Similar to the proof of  (6) in Theorem 2.5 we obtain
 $$||D(x_{\Lambda})|| \geq \pi n, \ n\geq 1.$$
The  last inequality contradicts the choice of $\pi\neq 0.$
Therefore  $\delta\equiv 0,$ i.e. from (11) we obtain that
$D_{n}=D_{a_{n}}.$ Since  $D_{a_{n}}$ is a spatial derivation and
the center of the algebra  $Z(e_nMe_n)$ coincides with $e_nZ,$ it
follows that $D_n$ is  $e_nZ$-linear. Thus
\begin{equation}
D_{n}(\lambda e_npe_n)=\lambda e_nD_{n}(e_npe_n)
\end{equation}
for all  $\lambda\in Z.$ Since the projection  $e_n$ is in $S_0(M,
\tau)$ and it commutes with $p$ we have  $$D_n(e_n pe_n)=D_n(e_n
p)=e_nD(e_n p)e_n=e_nD(e_n)pe_n+e_nD(p)e_n=$$$$=e_nD(e_n)e_n
p+e_nD(p)e_n=e_nD(p)e_n,$$ i.e.
\begin{equation}
\lambda D_n(e_n pe_n)=\lambda e_nD(p)e_n.
\end{equation}
In a similar way we obtain

 \begin{equation}
 D_n(\lambda e_n pe_n)=e_nD(\lambda p)e_n.
\end{equation}
Now (14), (15) and  (16) imply
$$e_{n}D(\lambda  p)e_{n}=e_{n}\lambda D(p)e_{n}$$ for all $n\in
\mathbb{N}$.

Set  $y=D(\lambda p)-\lambda D(p).$ Then  $e_n ye_n=0.$ From
$e_1=p_1=p,$  we have   $pyp=0.$ By Lemma 5.7 we have
$p^{\perp}yp^{\perp}=0.$ Multiplying the equality  $e_n ye_n=0$ by
 $p$ from the left side we obtain  $pye_n=0$ for all $n\in\mathbb{N}.$
 Since   $e_n\uparrow p_0+p,$ it follows that
 $py(p_0+p)=0,$ i.e.  $pyp_0=0.$ Since  $p_0$ is an
 arbitrary projection with the central cover $\textbf{1}$
such that $p_0\leq p^{\perp}$ and such that the von Neumann
algebra $p_0Mp_0$ is of type  I$_{\aleph_{0}},$  we obtain that
$pyp^{\perp}=0.$

 Similarly  $p^{\perp}yp=0.$ Therefore
 $$pyp=pyp^{\perp}=p^{\perp}yp=p^{\perp}yp^{\perp}=0$$ and hence
 $$y=pyp+pyp^{\perp}+p^{\perp}yp+p^{\perp}yp^{\perp},$$ i.e. $D(\lambda
p)=\lambda D(p).$ The proof is complete. $\blacksquare$

\textbf{Lemma  5.9.} \emph{Suppose that   $\lambda\in Z$ and
$x\in S_0(M, \tau).$ Then  }
 $$D(\lambda x)=\lambda D(x).$$

Proof. Case (i).  $x=p$ is a projection and
\begin{equation}
p=\sum\limits_{i=1}^{k}p_i,
\end{equation}
 where  $p_i,\, i=\overline{1, k}$
are mutually orthogonal abelian projections with
$\tau(p_i)<\infty.$ By Lemma  5.8 we have  $D(\lambda
p_i)=\lambda D(p_i).$ Therefore
$$D(\lambda p)=D(\lambda\sum\limits_{i=1}^{k}p_i)=\sum\limits_{i=1}^{k}D(\lambda
p_i)=\sum\limits_{i=1}^{k}\lambda D(p_i)=\lambda
D(\sum\limits_{i=1}^{k}p_i)=\lambda D(p),$$ i.e.
$$D(\lambda p)=\lambda D(p).$$

Case (ii).  $x=p$ is a projection with $\tau(p)<\infty.$ Then
$pMp$ is a finite von Neumann algebra of type I, and therefore
there exists a sequence of mutually orthogonal central projections
 $\{z_n\}$ such that  each $p_n=z_n p$ is a projection of the form (17). From the above case we have
 $D(\lambda p_n)=\lambda D(p_n).$ This and Lemma  5.6 imply that
 $$z_n D(\lambda p)=D(\lambda z_n p)=D(\lambda p_n)=\lambda D(p_n)=\lambda D(z_n p)=\lambda z_n D(p).$$
i.e.
 $$z_n D(\lambda p)=z_n\lambda  D(p)$$ for all $n.$
Therefore
$$D(\lambda p)=\lambda  D(p)$$

Case (iii). Let  $x\in S_0(M, \tau)$ be an element such that
$xp=x$ for some projection $p$ with $\tau(p)<\infty.$ Then
$$D(\lambda x)=D(\lambda x p)=D(x \lambda p)=D(x)\lambda p+xD(\lambda
p)=$$$$=D(x)\lambda p+x\lambda D(p)=\lambda(D(x)p+xD(p))=\lambda
D(xp)=\lambda D(x),$$ i.e. $D(\lambda x)=\lambda D(x).$

Case (iv).   $x$ is an arbitrary element from $S_0(M, \tau).$
Take a projection $p$ with  finite trace $\tau(p).$ Put $x_0=xp.$ From
the case (iii) we have  $D(\lambda x_0)=\lambda D(x_0).$ Now one
has
$$D(\lambda x_0)=D(\lambda x p)=D(\lambda x)p+\lambda xD(p),$$
i.e.
$$D(\lambda x)p=D(\lambda x_0)-\lambda xD(p).$$
On the other hand
$$D(\lambda x_0)=\lambda D(x_0)=\lambda D(xp)=\lambda D(x)p+\lambda x D(p),$$
i.e.
$$\lambda D(x)p=D(\lambda x_0)-\lambda x D(p).$$
Therefore   $\lambda D(x)p=D(\lambda x)p.$ Since  $p$ is an
arbitrary with $\tau(p)<\infty,$ this implies
$$D(\lambda x)=\lambda D(x).$$ The proof is complete. $\blacksquare$

Proof of Theorem 5.5.

By Lemma  5.9 the derivation  $D: S_0(M, \tau)\rightarrow S_0(M,
\tau)$ is  $Z$-linear. By Theorem  5.1 $D$ is spatial and moreover
$$D(x)=ax-xa,\,x\in S_0(M, \tau)$$ for an appropriate  $a\in S(M,
\tau).$ The proof is complete. $\blacksquare$

Now we can describe the structure of derivations on the
algebra $S_0(M, \tau)$ of $\tau$-compact operators with respect to
a type I von Neumann algebra $M$ with a faithful normal
semi-finite trace $\tau.$

Let $M$ be a type  I von Neumann algebra and   let $z_0$ be the  central
projection in $M$ such that $z_0M$ is a finite von Neumann algebra and
$z_0^{\bot}M$ is a von Neumann algebra of type  I$_{\infty}.$
Consider a derivation  $D$ on  $S_0(M, \tau)$ and let $\delta$
be its restriction onto the center  $Z(S_0(M, \tau)).$ By Proposition 1.2
we have $z_0^{\bot}Z(S_0(M, \tau))=\{0\},$ and therefore
$z_0^{\bot}\delta\equiv 0,$ i.e. $\delta=z_0\delta.$

By Lemma 4.2 $\tau(z_\delta)<\infty$ and therefore the derivation $D_\delta$  defined
 in (2) maps   $z_0S_0(M, \tau)$ into itself. Consider its extension $D_\delta$ on $S_0(M,
\tau)=z_0S_0(M, \tau)\oplus z_0^{\bot}S_0(M, \tau)$ which is
defined as
\begin{equation}
D_\delta(x_1+x_2):=D_\delta(x_1),\,x_1\in z_0S_0(M, \tau),x_2\in z_0^{\bot}S_0(M,
\tau).
\end{equation}
Similar to the cases of the algebras $LS(M),$ $S(M)$ and $S(M, \tau)$
 for a finite von Neumann algebra $M$ of type I, every derivation on the algebra
  $S_0(M, \tau)$ admits the decomposition $D=D_a+D_\delta.$

The following is the main result of this section, which gives
the general form of derivations on the algebra $S_0(M, \tau)$ (cf. \cite{Alb4}).

\textbf{Theorem  5.10.} \emph{Let  $M$ be a type  $I$ von Neumann
algebra with a faithful normal semi-finite trace $\tau.$ Each
derivation  $D$ on  $S_0(M, \tau)$ can be uniquely represented in
the form
\begin{equation}
D=D_{a}+D_{\delta,}
\end{equation}
where  $D_{a}$ is a spatial derivation implemented by an element
$a\in S(M, \tau),$ and $D_{\delta} $ is a derivation of the form
(18), generated by a derivation $\delta$ on the center of $S_0(M,
\tau)$.}

Similar to Corollary 4.5 we obtain the following result.

 \textbf{Corollary  5.11.} \emph{Under the conditions of Theorem 5.10 a derivation
  $D$ on  the algebra $S_0(M, \tau)$ is spatial if and only if $D$ is continuous in the measure topology
  $t_\tau.$}

Finally from Theorems 2.7, 3.6, 4.4, 5.10 and from \cite[Theorem 3.4]{Ber}
we obtain the following corollary

\textbf{Corollary 5.12.} \emph{Let $M$ be a type I von Neumann algebra. The following
conditions are equivalent:}

\emph{(i) Every derivation on the algebra $LS(M)$ (resp. $S(M),$ $S(M, \tau)$)
is inner.}

\emph{(ii) Every derivation on the algebra  $S_0(M, \tau)$
is spatial.}

\emph{(iii) The center of the type I$_{fin}$ part of $M$ is atomic.}

\begin{center} \textbf{6.  An application to the description of the first cohomology
group}\end{center}

Let  $A$ be an algebra. Denote by  $Der(A)$ the space of all
derivations (in fact it is a Lie algebra with respect to the
commutator), and denote  by $InDer(A)$ the subspace of all inner
derivations on $A$ (it is a Lie ideal in $Der(A)$).

 The factor-space   $H^{1}(A)=Der(A)/InDer(A)$ is called the first (Hochschild) cohomology group of the algebra
  $A$ (see \cite{Dal}). It is clear that  $H^{1}(A)$ measures how much the space of all derivations
  on $A$ differs from the space on inner derivations.

The  following result shows that the first cohomology groups of the algebras $LS(M),$ $S(M)$ and $S(M, \tau)$
are completely determined by the corresponding cohomology groups of their centers (cf. \cite[Corollary 3.1]{Ber}).

\textbf{Theorem 6.1.} \emph{Let  $M$ be a type I von Neumann
algebra with the center $Z$ and a faithful normal semi-finite
trace $\tau.$ Suppose that  $z_0$ is a central projection such
that  $z_0M$ is a finite von Neumann algebra, and  $z_0^{\bot}M$
is of type  I$_{\infty}.$ Then}

\emph{a)  $H^{1}(LS(M))=H^{1}(S(M))\cong H^{1}(S(z_0Z));$}

\emph{b) $H^{1}(S(M, \tau))\cong H^{1}(S(z_0Z, \tau_0)),$
where $\tau_0$ is the restriction of  $\tau$ onto
$z_0Z.$ }

Proof. It immediately follows from Theorems 2.7, 3.6 and 4.4. $\blacksquare$

Further we need the following property of the algebra of
$\tau$-compact operators from  \cite{Str}:
\begin{equation}
S(M, \tau)=M+S_0(M, \tau).
\end{equation}
Set $C(M, \tau)=M\cap S_0(M, \tau)$ and consider  $M/(C(M,
\tau)+Z)$ -- the factor space of  $M$ with respect to the space
$C(M, \tau)+Z.$

For  $D_1, D_2\in Der(S_0(M, \tau))$ put
$$D_1 \sim D_2 \Leftrightarrow D_1-D_2\in InDer(S_0(M, \tau)).$$
Suppose that   $D_1 \sim D_2.$ From Theorem 5.10 these derivation can be
represented in the form (19):
$$D_1=D_a+D_\delta,\, D_2=D_b+D_\sigma.$$
Since $D_1-D_2=D_c,$ where $c\in S_0(M, \tau)\subset S(M, \tau),$
from the uniqueness of a the representation in the form (19) it follows
that  $D_a-D_b\in InDer(S_0(M, \tau))$ and $D_\delta=D_\sigma.$
Therefore  $\delta\equiv\sigma$ and
\begin{equation}
a-b\in S_0(M, \tau)+Z(S(M, \tau)).
\end{equation}
 According to  (20) we have
$$a=a_1+a_2, a_1\in M,\,a_2\in S_0(M, \tau),$$
$$b=b_1+b_2, b_1\in M,\,b_2\in S_0(M, \tau).$$
From  (21) it follows that
$$a_1-b_1\in (b_2-a_2)+S_0(M, \tau)+Z(S(M,
\tau))\subset S_0(M, \tau)+Z(S(M, \tau)).$$ Since  $a_1, b_1\in
M,$ we have that   $$a_1-b_1\in (S_0(M, \tau)+Z(S(M, \tau)))\cap
M\subset C(M, \tau)+Z$$ because $Z(S(M, \tau))\cap M=Z$ (cf. Proposition 1.2). Therefore
$$D_1 \sim D_2 \Leftrightarrow a_1-b_1\in C(M, \tau)+Z,\, \delta\equiv\sigma.$$

Thus we have the following result.

\textbf{Theorem 6.2.} \emph{Let  $M$ be a type I von Neumann
algebra with the center $Z$ and a faithful normal semi-finite
trace $\tau.$ Suppose that  $z_0$ is a central projection such
that  $z_0M$ is a finite von Neumann algebra, and  $z_0^{\bot}M$
is of type  I$_{\infty}.$ Then the group  $H^{1}(S_0(M, \tau))$ is
isomorphic with the group  $M/(C(M, \tau)+Z)\oplus H^{1}(S_0(z_0Z,
\tau_0)),$ where $\tau_0$ is the restriction of  $\tau$ onto
$z_0Z.$ In particular, if  $M$ is of type  I$_{\infty},$ then $H^{1}(S_0(M,
\tau))\cong M/(C(M, \tau)+Z).$ }

\textbf{Remark 3.} In the algebras $S(M, \tau)$ and $S_0(M, \tau)$ equipped with the measure topology $t_\tau$ one can consider another possible cohomology theories. Similar to \cite{KR} consider the space $Der_c(A)$ of all continuous derivation on a topological algebra $A$ and define the first cohomology group $H_{c}^{1}(A)=Der_c(A)/InDer(A).$

Under these notations the above results and Corollaries 4.6 and 5.11 imply the following result (cf. \cite[Theorem 4.4]{KR}).

\textbf{Corollary 6.3.} \emph{Let  $M$ be a type I von Neumann
algebra with the center $Z$ and a faithful normal semi-finite
trace $\tau.$ Consider the topological algebras $S(M, \tau)$ and $S_0(M, \tau)$ equipped with the measure topology. Then $H^{1}_c(S(M, \tau))=\{0\}$ and $H^{1}_c(S_0(M, \tau))\cong M/(C(M, \tau)+Z).$ }

\vspace{1cm}

\textbf{Acknowledgments.} \emph{The second and third named authors
would like to acknowledge the hospitality of the $\,$ "Institut
f\"{u}r Angewandte Mathematik",$\,$ Universit\"{a}t Bonn
(Germany). This work is supported in part by the DFG 436 USB
113/10/0-1 project (Germany) and the Fundamental Research
Foundation of the Uzbekistan Academy of Sciences.}

\end{document}